\documentclass[hidelinks,onefignum,onetabnum]{siamart250211}
\usepackage{lipsum}
\usepackage{amsfonts}
\usepackage{url}
\usepackage{mathtools}
\usepackage{cleveref}
\usepackage{doi}
\usepackage{tikz}
\usepackage{pgfplots}
\pgfplotsset{compat=1.18}
\usetikzlibrary{plotmarks}
\usepackage{thm-restate}
\usepackage{graphicx}
\usepackage{epstopdf}
\usepackage{algorithmic}
\ifpdf
  \DeclareGraphicsExtensions{.eps,.pdf,.png,.jpg}
\else
  \DeclareGraphicsExtensions{.eps}
\fi

\newsiamremark{remark}{Remark}
\newsiamthm{assumption}{Assumption}

\headers{Posterior error bounds for prior-driven balancing}{J. König and H. C. Lie}

\title{Posterior error bounds for prior-driven balancing in linear Gaussian inverse problems}

\author{Josie König\thanks{Institut für Mathematik, Universität Potsdam, Potsdam OT Golm 14476, Germany (\email{josie.koenig@uni-potsdam.de},  \email{han.lie@uni-potsdam.de}).}  
\and Han Cheng Lie\footnotemark[2]}

\usepackage{amsopn}

\newcommand{\pr}[0]{\mathrm{pr}}
\newcommand{\pos}[0]{\mathrm{pos}}
\newcommand{\obs}[0]{\mathrm{obs}}
\newcommand{\unknown}{\mbf{p}}
\newcommand{\data}{\mbf{m}}
\newcommand{\noise}{\pmb{\epsilon}}
\newcommand{\fwdmap}[0]{\mbf{G}}
\newcommand{\hatFwdmap}[0]{\widehat{\mbf{G}}}
\newcommand{\PrMn}{\pmb{\mu}_{\pr}}

\newcommand{\PrCov}{\mathbf\Gamma_{\pr}}

\newcommand{\ObsCov}{\mathbf\Gamma_{\obs}}
\newcommand{\PosMn}{\pmb{\mu}_{\pos}}
\newcommand{\hatPosMn}{\widehat{\pmb{\mu}}_{\pos}}
\newcommand{\PosCov}{\mathbf\Gamma_{\pos}}
\newcommand{\hatPosCov}{\widehat{\mathbf\Gamma}_{\pos}}

\newcommand{\R}[0]{\mathbb{R}}
\newcommand{\N}[0]{\mathbb{N}}
\newcommand{\mbf}[1]{\mathbf{#1}}
\newcommand{\norm}[1]{\ensuremath{\lVert #1 \rVert}}
\newcommand{\Norm}[1]{\ensuremath{\left\lVert #1 \right\rVert}}

\DeclareMathOperator*{\range}{ran} 

\ifpdf
\hypersetup{
  pdftitle={Posterior error bounds for prior-driven balancing},
  pdfauthor={J. König and H. C. Lie}
}
\fi

\begin{document}
\maketitle

\begin{abstract}
In large-scale Bayesian inverse problems, it is often necessary to apply approximate forward models to reduce the cost of forward model evaluations, while controlling approximation quality. In the context of Bayesian inverse problems with linear forward models, Gaussian priors, and Gaussian noise, we use perturbation theory for inverses to bound the error in the approximate posterior mean and posterior covariance resulting from a linear approximate forward model. We then focus on the smoothing problem of inferring the initial condition of linear time-invariant dynamical systems, using finitely many partial state observations.
For such problems, and for a specific model order reduction method based on balanced truncation, we show that the impulse response of a certain prior-driven system is closely related to the prior-preconditioned Hessian of the inverse problem.
This reveals a novel connection between systems theory and inverse problems.
We exploit this connection to prove the first a priori error bounds for system-theoretic model order reduction methods applied to smoothing problems.
The bounds control the approximation error of the posterior mean and covariance in terms of the truncated Hankel singular values of the underlying system.
\end{abstract}



\section{Introduction}

In a linear statistical inverse problem, one assumes an `observation model' that describes the relationship $\data = \fwdmap(\unknown)+\noise$ between parameters $\unknown$ and noisy observations $\data$.
The parameter-to-observable map or `forward model' $\fwdmap$ is assumed to be continuous and linear, and $\noise$ is a random variable. Both $\fwdmap$ and the distribution of $\noise$ are known.
The goal is to infer, for a given realization of $\data$, the data-generating parameter $\unknown$.
Under the Bayesian approach, one models the unknown $\unknown$ as a random variable with a prior distribution.
We consider Gaussian linear inverse problems with finite-dimensional parameters and data, i.e. linear statistical inverse problems where the prior distribution of the unknown $\unknown$ and of the noise $\noise$ are finite-dimensional, Gaussian, and mutually independent; see e.g. \cite[Sections 2.1-2.2]{Flath2011}.

Linear Gaussian inverse problems appear in many scientific and engineering disciplines, e.g. in the analysis step of the Kalman filter in the context of data assimilation, and also in classical parameter estimation problems, such as inferring the initial condition of a linear dynamical system with constant coefficients. 
A key advantage of linear Gaussian inverse problems is that their solutions---i.e. the posterior distributions---are Gaussian measures, and hence completely determined by their mean and covariance.
The posterior mean vector and covariance matrix can be computed exactly, using the formulas for the conditional mean and conditional covariance of a multivariate Gaussian.
For example, the posterior covariance is the Schur complement of the marginal covariance of the data.

In practice, the numerical solution of inverse problems requires approximating the exact forward map $\mbf{G}$ to mitigate computational costs in large-scale settings, even for linear forward models. In high-dimensional parameter spaces, naive approximations of forward models may be too costly to apply repeatedly, and additional steps are needed to obtain practical algorithms, see e.g. \cite{Chung2017}.
One often applies a combination of dimension reduction and model order reduction (MOR), in order to obtain lower-dimensional or reduced forward models that are practical for computational inference, for example as in \cite{Flath2011}.

The use of approximate forward models naturally leads to the question of measuring the quality of the corresponding approximate posterior, when the prior and noise are fixed. 
For linear Gaussian inverse problems, changing the linear forward model yields an approximate posterior that is Gaussian, and thus it suffices to quantify the approximation quality of the posterior mean and covariance. 
If the approximate forward map $\hatFwdmap$ is obtained via an MOR technique, then it is desirable to control the error in the approximate posterior in terms of objects that are key descriptors of the MOR technique.

\subsection{Contributions}
\label{subsection_contributions} 
We use perturbation theory for inverses \cite{Wedin1973} to prove the first a priori error bounds on the posterior covariance matrix and on the posterior mean vector associated to a linear perturbation of the forward model, for linear Gaussian inverse problems. 
An important aspect of these bounds is that they allow for singular prior covariances.
This offers important flexibility for large-scale inference, by allowing to focus the prior-to-posterior update on proper subspaces of high-dimensional parameter spaces, e.g. subspaces spanned by an ensemble of particles.
Another important aspect of these bounds is that the error is controlled in terms of the corresponding error in square roots of the prior-preconditioned Hessian matrix, and thus account for regularization by both the prior and the noise.

Second, we reveal a previously unknown connection between systems theory and the theory of linear Gaussian inverse problems, in the context of the smoothing problem from data assimilation for linear time-invariant (LTI) dynamical systems.
Namely, we show that the impulse response functions associated to discrete-time measurements and prior-driven balanced truncation form the rows of the prior-preconditioned Hessian.

Third, we exploit the above-mentioned connection to prove the first a priori error bounds for the approximate posterior mean and posterior covariance associated to a class of system-theoretic MOR methods for data assimilation based on balanced truncation. 
Specifically, for the balancing methods applied to prior-driven systems described in \cite{2025PD-BT}, we prove bounds for standard and time-limited balanced truncation separately.
The bounds are given in terms of the truncated Hankel singular values of the prior-driven system.
These singular values are system-theoretic invariants that provide insight into selecting the truncation rank for the approximation. 

Numerical experiments illustrate our theoretical results and provide the first implementation of time-limited balanced truncation for prior-driven systems.

\subsection{Related literature}
\label{subsection_related_literature}

Some of the earliest approximations of the forward model for Bayesian inverse problems result from discretizations of spatial or spatiotemporal fields arising from PDEs \cite{Cotter2010BIPforPDEs,Dashti2011UQEllipticIVPs}.
Later approximations take into account computational cost, e.g. methods based on generalized polynomial chaos \cite{Yan2017ConvSurrogateBIPs} or the reduced basis method \cite{Chen2015LinSGRBBIPs}.
In \cite{Flath2011}, a method for approximating the Gaussian posterior using a low-rank approximation of the prior-preconditioned Hessian was proposed.
Low-rank matrix approximations for Gaussian approximations of posteriors were also proposed in \cite{Bonnabel2024}, but in the context of matrix differential equations.
The posterior covariance field of Gaussian processes is analyzed in \cite{Cai2025}, but with a focus on spatial features such as the kernel bandwidth and the location of observations. There is exhaustive literature on low-rank approximations for general kernel matrices, see e.g. \cite{Cai2023Data-DrivenLRApprox,Heilig2022KernelApprox,Kumar2009SamplingSpectral,Zhang2010ClusteredNystrom}, but low-rank approximations for general kernel matrices are not the focus of this work. Instead, we focus on approximations of the forward model and the resulting error in the Bayesian posterior derived from these approximations.

The most recent error bounds for Bayesian inverse problems appear to be those in \cite{CvetkovicLie2025,Sprungk2020LocalLipschitz}; see in particular \cite[Remark 4]{Sprungk2020LocalLipschitz} for comments about perturbations of the forward model.
Because these bounds do not assume linear forward models or Gaussian priors, the bounds are given in terms of distances or divergences of the approximate posterior measure, and thus are distinct from the bounds we prove below.
In \cite{SanzAlonso2024}, the focus is on error bounds for finite element- or graph-based discretizations of a linear forward model, instead of arbitrary linear perturbations of the forward model.
In \cite{Spantini2015Optimal}, the low-rank approximation from \cite{Flath2011} was shown to yield information-theoretically optimal approximations of the covariance and the mean separately; the approximation error is an explicit function of the truncated eigenvalues of the prior-preconditioned Hessian.
In \cite{CarereLie2025a,CarereLie2025b}, these optimality properties were shown to be intrinsic properties of linear Gaussian inverse problems on Hilbert spaces, in the sense of not depending on the discretization method used.
However, these optimal low-rank approximations are obtained by pre-composing the forward model with a projector, and not by additive linear perturbations of the forward model.

Research into system-theoretic MOR methods for linear Gaussian inverse problems was initiated only very recently, with the focus being on the smoothing problem in the context of data assimilation for high-dimensional LTI systems. 
To render such inverse problems computationally tractable, \cite{Qian2021Balancing} proposed a system-theoretic MOR method that aims at recovering the optimal low-rank approximation (OLR) proposed in \cite{Flath2011} for linear Gaussian inverse problems. 
A consistency argument in the limit of infinite-time horizon and infinitely frequent observations is given for the posterior, but no error bound is given.
In \cite{Koe2022TLBTDA}, the system-theoretic MOR method from \cite{Qian2021Balancing} was applied to linear Gaussian inverse problems and to the 4D-Var method from data assimilation, and a time-limited balanced truncation method was proposed. The time-limited method was further analyzed and compared to $\alpha$-bounded balanced truncation for unstable discrete systems in \cite{koenig2022TLBT4DVAR}.
In \cite{2025PD-BT} a MOR method based on balanced truncation (BT) of the auxiliary prior-driven system (PD-BT) is proposed, and an error bound on the output approximation is given.
For the methods described above, no error bounds for the posterior mean and covariance are available in the literature.

\subsection{Outline}

We introduce the key notation for this paper in \cref{section_notation}. \Cref{sec:Preliminaries} describes the setting and the smoothing problem for LTI systems.
\Cref{sec:Lipschitz bounds} discusses the local Lipschitz stability of Gaussian posteriors, and presents error bounds for the posterior mean and covariance in terms of the prior-preconditioned Hessian in \cref{theorem03}.
In \Cref{sec:PD_BT_bound}, we describe a novel connection between systems theory and the theory of linear Gaussian inverse problems, and use this connection to prove \cref{thm:PD-BT bounds} and \cref{thm:PD-TLBT bounds}, which yield the first a priori error bounds on posterior approximations for prior-driven balanced truncation (PD-BT).
We show and discuss numerical results that illustrate these bounds in \Cref{sec:Numerics}.
After discussing the possibility for proving analogous bounds for other system-theoretic MOR methods in \cref{sec:MOR Alternatives}, we conclude in \Cref{sec:Conclusion}.
Auxiliary results and the proof of \cref{theorem03} are given in \cref{section_proofs_auxiliary_lemmas} and \cref{section_proofs_local_lipschitz_bounds}, respectively.

\subsection{Notation}
\label{section_notation}
For $\mbf{A}\in\R^{m\times n}$, we denote the range or column span of $\mbf{A}$ by $\range \mbf{A}$, the Moore--Penrose inverse by $\mbf{A}^\dagger\in\R^{n\times m}$, and the $p$-Schatten norm of $\mbf{A}$ for $1\leq p \leq \infty$ by $\norm{\mbf{A}}_p$.
In particular, $\norm{\mbf{A}}_F\equiv \norm{\mbf{A}}_2$ denotes the Frobenius or Hilbert--Schmidt norm, while $\norm{\mbf{A}}_\infty$ denotes the spectral norm $\norm{\mbf{A}}_\infty=\sup\{ \norm{\mbf{A}\mbf{x}}_2:\norm{\mbf{x}}_2=1\}$.
For a symmetric matrix $\mbf{A}\in\R^{d\times d}$, we shall write $\mbf{A}\succeq 0$ and $\mbf{A}\succ 0$ if $\mbf{A}$ is positive semidefinite and positive definite, respectively.
If $\mbf{A}\succeq 0$, then there exists a unique $\mbf{B}\succeq 0$ that satisfies $\mbf A=\mbf B \mbf B$; see e.g. \cite[p. 405, Theorem 7.2.6]{HornJohnson1985}.
We shall denote this `principal non-negative square root' of $\mbf A$ by $\mbf A^{1/2}$.
For a vector $\mbf{x}\in\R^{m}$, $\norm{\mbf{x}}_p$ denotes the $\ell_p$ norm.

\section{Preliminaries}\label{sec:Preliminaries}
We introduce the general setting and solution of linear Gaussian inverse problems in \cref{subsec:problem setting}, before describing the special case where the forward model arises from the smoothing problem for a linear time-invariant (LTI) system in \cref{subsec:LTI FOM}. Prior-driven balancing, a system-theoretic MOR method for such systems, is introduced in \cref{subsec:PD-BT}.
\subsection{Setting}\label{subsec:problem setting}
\begin{assumption}
\label{assumption_main}
We assume the observation model
    \begin{align}\label{eq:generic_IVP}
    \data = \fwdmap\unknown+\noise, 
\end{align}
for a forward model $\fwdmap \in \R^{d_\obs \times d}$ for some $d, d_\obs\in \mathbb{N}$, and centered $\R^{d_\obs}$-valued Gaussian observation noise $\noise\sim \mathcal{N}(\mbf{0},\ObsCov)$, where $\ObsCov\succ 0$.
We place a Gaussian prior on the $\R^{d}$-valued unknown $\unknown\sim \mathcal{N}(\PrMn,\PrCov)$, where $\PrCov\succeq 0$.
\end{assumption}
\Cref{assumption_main} implies that $\ObsCov$ is invertible, and allows for $\PrCov$ to be singular.

Under \Cref{assumption_main}, the solution to the Bayesian inverse problem of inferring the parameter $\unknown$ from the realization of the data $\data$ is a Gaussian measure $\unknown|\data\sim\mathcal{N}(\PosMn(\data),\PosCov)$ with
\begin{subequations}
\label{eq:Posterior}
    \begin{align}
    \PosMn(\data) &= \PrMn + \PosCov\fwdmap^\top\ObsCov^{-1}(\data-\fwdmap\PrMn) \in \R^{d},
    \label{eq:pos_mean}
    \\
    \PosCov &= \PrCov - \PrCov \fwdmap^\top(\ObsCov+\fwdmap\PrCov \fwdmap^\top)^{-1}\fwdmap\PrCov\in \R^{d \times d},
    \label{eq:pos_covariance}
\end{align}
\end{subequations}
see e.g. \cite[equations (6.14), (6.15)]{Stuart2010}.
The Hessian of the negative log-likelihood with respect to the unknown parameter variable $\unknown$ is $\mbf{H}=\fwdmap^\top\ObsCov^{-1}\fwdmap\in \R^{d \times d}$. Let $\mbf{L}_\pr\in \mathbb{R}^{d\times s}$ be a possibly non-symmetric square root of $\PrCov$, i.e. $\PrCov=\mbf{L}_\pr\mbf{L}_\pr^\top$. We shall then refer to the operator $\mbf{L}_\pr^\top\mbf{H}\mbf{L}_\pr\in\mathbb{R}^{s\times s}$ as the \textit{prior-preconditioned Hessian} and $\ObsCov^{-1/2}\fwdmap\mbf{L}_\pr\in \mathbb{R}^{d_{\rm obs}\times s}$ as its corresponding \textit{square root}.

If instead of $\fwdmap$ we use an approximation $\hatFwdmap$ of $\fwdmap$ and leave the data $\data$, prior $\mathcal{N}(\PrMn,\PrCov)$ and observation noise $\mathcal{N}(\mbf{0},\ObsCov)$ unchanged, then this results in an approximation $\mathcal{N}(\hatPosMn(\data),\hatPosCov)$ of the `exact' Gaussian posterior $\mathcal{N}(\PosMn(\data),\PosCov)$, where the approximate posterior mean $\hatPosMn$ and approximate covariance $\hatPosCov$ are obtained by replacing $\fwdmap$ in \eqref{eq:Posterior} with $\hatFwdmap$:
\begin{subequations}
\label{eq:PD-BT PosBoth}
\begin{align}
    \hatPosMn(\data) &=\PrMn + \widehat{\mbf{\Gamma}}_{\pos}\widehat{\fwdmap}^\top\ObsCov^{-1}(\data-\widehat{\fwdmap}\PrMn),
    \label{eq:pos_mean_approx}
    \\
    \hatPosCov &= \PrCov - \PrCov \widehat{\fwdmap}^\top(\ObsCov+\widehat{\fwdmap}\PrCov \widehat{\fwdmap}^\top)^{-1}\widehat{\fwdmap}\PrCov.
    \label{eq:pos_covariance_approx}
\end{align}
\end{subequations}

\subsection{Forward model from an LTI system}\label{subsec:LTI FOM}
We will consider a special case of the linear Gaussian inverse problem where the initial condition of a high-dimensional linear time-invariant (LTI) system is the unknown. 
Let the matrices $\mbf{A} \in \R^{d \times d}$ and $\mbf{C} \in \R^{d_{\mathrm{out}} \times d}$, the state $\mbf{x}(t)\in\R^{d}$ and the output $\mbf{y}(t)\in\R^{d_{\rm out}}$ satisfy
\begin{subequations}\label{eq:BayesianSystem}
\begin{align}
\begin{split}\label{eq:BayesianSystem_nonoise}
    \dot{\mbf{x}}(t)  &= \mbf{A}\mbf{x}(t),\enspace\quad \mbf{x}(0)=\unknown\sim \mathcal{N}(\mbf{0},\PrCov),\\
    \mbf{y}(t) &= \mbf{C}\mbf{x}(t).
\end{split}
\end{align}
The goal is to infer $\unknown$ from $n$ output measurements $\data = \begin{bmatrix} \data_1^\top, \ldots, \data_n^\top \end{bmatrix}^\top\in \R^{nd_{\mathrm{out}}}$ taken at discrete observation times $t_{1} < \ldots < t_{n}$ and given by
\begin{align}\label{eq:BayesianSystem_noise}
    \data_k &= \mbf{y}(t_k) + \noise_k= \mbf{C}\mbf{x}(t_k) + \noise_k, \enspace \noise_k \sim \mathcal{N}(\mbf{0},\mbf{\Gamma}_\epsilon), \enspace k=1,\ldots,n,
\end{align}
where the $(\noise_k)_k$ are independent and identically distributed, and $\mbf{\Gamma}_\epsilon\in\R^{d_{\rm out}\times d_{\rm out}}\succ 0$.
\end{subequations}

The problem of inferring the unknown initial condition $\mathbf{p}$ from  a series of measurements taken after the initial time may also be referred to as the \textit{Bayesian smoothing problem}. 
It is a linear Gaussian inverse problem~\eqref{eq:generic_IVP} with $\PrMn=\mbf{0}$ and linear forward model and observation covariance matrix as follows:
\begin{align}\label{eq:BLIPproperties}
\fwdmap= \begin{bmatrix}\mbf{C}e^{\mbf{A}t_1}\\\vdots\\\mbf{C}e^{\mbf{A}t_n}\end{bmatrix}\in \R^{d_{\rm obs}\times d}, \quad\text{and}\quad \ObsCov= \begin{bmatrix} \mbf{\Gamma}_\epsilon\\ & \ddots \\ && \mbf{\Gamma}_\epsilon \end{bmatrix}\in \R^{d_{\rm obs}\times d_{\rm obs}},
\end{align}
where $d_{\rm obs} = n\cdot d_{\rm out}$. The posterior is given by~\eqref{eq:Posterior} with $\PrMn=\mbf{0}$ and thus $\pmb{\mu}_{\pos}(\data) =\PosCov\fwdmap^\top\ObsCov^{-1}\data$.

In large-scale scenarios with high state dimension $d$ and possibly high observation dimension $d_{\rm obs}$, directly constructing the operator $\fwdmap$ can be impractical. Instead, the relationship between the initial conditions and the measurements typically remains implicit via the time integration of the high-dimensional dynamical system described in~\eqref{eq:BayesianSystem_nonoise}. To alleviate the computational burden of evaluating $\fwdmap$ in such situations, MOR techniques are used to replace the high-dimensional model~\eqref{eq:BayesianSystem_nonoise} with a lower-dimensional \textit{reduced model}. In this work, we will focus on the prior-driven balancing approach that was introduced in \cite{2025PD-BT}.

\subsection{Prior-driven balancing for Bayesian inference}\label{subsec:PD-BT}
Prior-driven balanced truncation (PD-BT) \cite{2025PD-BT} considers the \textit{prior-driven} homogeneous LTI system,
\begin{align}\label{eq:PriorDrivenSys}
    \begin{split}
        \dot{\mbf{x}}(t)  &= \mbf{A}\mbf{x}(t)+\mbf{L}_\pr\mbf{u}(t), \enspace \mbf{x}(0) = \mbf{0},\\
        \mbf{y}_\epsilon(t) &= \mbf{\Gamma}_\epsilon^{-1/2}\mbf{C}\mbf{x}(t),
    \end{split}
\end{align} 
with impulse response
\begin{align}\label{eq:FullImpulseResponse}
    \mbf{h}(t)= \mbf{\Gamma}_\epsilon^{-1/2}\mbf{C}e^{\mbf{A}t}\mbf{L}_\pr.
\end{align}
The impulse response of a system is a key system-theoretic property that completely defines the behavior of an LTI system \cite[Chapter 4.1]{Antoulas2005Book}. The $L_2$-norm of the impulse response is defined as
\begin{align}\label{eq:L2Norm}
    \norm{\mbf{h}}_{L_2} = \sqrt{\int_0^\infty\mathrm{trace}\left(\mbf{h}^\top (t)\mbf{h}(t)\right)\,\text{d}t} =  \sqrt{\int_0^\infty\norm{\mbf{h}(t)}_F^2\,\text{d}t},
\end{align}
and the time-limited $L_2$-norm at time $\mathcal{T}\in \R_{>0}$ as
\begin{align}\label{eq:TLL2Norm}
    \norm{\mbf{h}}_{L_{2,\mathcal{T}}} = \sqrt{\int_0^\mathcal{T}\mathrm{trace}\left(\mbf{h}^\top (t)\mbf{h}(t)\right)\,\text{d}t} =  \sqrt{\int_0^\mathcal{T}\norm{\mbf{h}(t)}_F^2\,\text{d}t}.
\end{align}

The idea of balanced truncation is to identify states that are simultaneously hard to reach because they require high-energy inputs and hard to observe because their associated outputs dissipate little energy. These unimportant states are then truncated from the state space to reduce its dimension. The reachability and observability energies of states in an LTI system are associated with two matrices, respectively: the reachability Gramian $\mathbf{P}$ and the observability Gramian $\mathbf{Q}$. For more information, we refer to \cite[Chapter 7]{Antoulas2005Book}.

We assume the system matrix $\mbf{A}$ to be stable, meaning that all eigenvalues of $\mbf{A}$ lie in the open left half-plane. Then, the system Gramians for PD-BT, i.e., \textit{standard} balanced truncation \cite{Moore1981BT,Mullis1976BT} of the system~\eqref{eq:PriorDrivenSys}, are given by 
\begin{align}\label{eq:Inf_Gramians}
  \mathbf{P}_\infty =\int_0^\infty e^{\mathbf{A} t}\PrCov e^{\mathbf{A}^\top t}\,\text{d}t, \quad\text{ and }\quad
  \mathbf{Q}_\infty =\int_0^\infty e^{\mathbf{A}^\top t}\mathbf{C}^\top\boldsymbol{\Gamma}_\epsilon^{-1}\mathbf{C}e^{\mathbf{A}t}\,\text{d}t.  
\end{align}
We can also consider the prior-driven analogue of \textit{time-limited} balanced truncation \cite{Wodek1990TLGramians}, which we shall refer to as `PD-TLBT' to distinguish it from PD-BT.
The system Gramians for PD-TLBT for an arbitrary, not necessarily stable system matrix $\mbf{A}$ are given by
\begin{align}\label{eq:TL_Gramians}
  &\mathbf{P}_\mathcal{T} =\int_0^\mathcal{T} e^{\mathbf{A} t}\PrCov e^{\mathbf{A}^\top t}\,\text{d}t, \quad\text{ and }\quad \mathbf{Q}_\mathcal{T} =\int_0^\mathcal{T} e^{\mathbf{A}^\top t}\mathbf{C}^\top\boldsymbol{\Gamma}_\epsilon^{-1}\mathbf{C}e^{\mathbf{A}t}\,\text{d}t.  
\end{align}
These Gramians are used for the computation of reduced bases $\mbf{W}_r,\mbf{V}_r\in\R^{d\times r}$. For this, let $\mathbf{P}=\mbf{P_*} = \mbf{LL}^\top$ and $ \mbf{Q}= \mbf{Q_*} = \mbf{RR}^\top$ for $*\in\{\mathcal T ,\infty\}$ be square root decompositions of the Gramians. Let $\mbf{R}^\top \mbf{L}=\mbf{U}\mbf{\Sigma}\mbf{Z}^\top$ be the singular value decomposition of the product $\mbf{R}^\top \mbf{L}$, and let $\mbf{U}_r\mbf{\Sigma}_r\mbf{Z}_r^\top$ be its rank-$r$ approximation. The diagonal elements of $\mbf{\Sigma}$ are called the \textit{Hankel singular values} (HSVs). The reduced bases are given by $\mbf{W}_r = \mbf{L}\mbf{Z}_r\mbf{\Sigma}^{-1/2}_r$ and $\mbf{V}_r =\mbf{R}\mbf{U}_r\mbf{\Sigma}^{-1/2}_r$ so that $\mbf{V}_r^\top\mbf{P}\mbf{V}_r=\mbf{\Sigma}_r=\mbf{W}^\top_r\mbf{Q}\mbf{W}_r$. The \textit{reduced} prior-driven system is then defined as
\begin{align}\label{eq:ReducedInhomBayesianSys}
    \begin{split}
        \dot{\mbf{x}}_r(t)  &= \mbf{A}_r\mbf{x}_r(t)+\mbf{L}_{\pr,r}\mbf{u}(t), \enspace \quad\mbf{x}(0) = \mbf{0},\\
         \mbf{y}_{\epsilon,r}(t) &= \mbf{\Gamma}_\epsilon^{-1/2}\mbf{C}_r\mbf{x}_r(t),
    \end{split}
    \end{align}
with $\mbf{x}_r\coloneqq \mbf{V}_r^\top\mbf{x}\in\mathbb{R}^r$, and reduced operators $\mbf{A}_r \coloneqq  \mbf{V}_r^\top\mbf{A}\mbf{W}_r\in\mathbb{R}^{r\times r}$, $\mbf{L}_{\pr,r} \coloneqq  \mbf{V}_r^\top \mbf{L}_\pr\in\mathbb{R}^{r\times s}$, and $\mbf{C}_r \coloneqq \mbf{C}\mbf{W}_r\in\mathbb{R}^{d_\mathrm{out}\times r}$. The reduced system~\eqref{eq:ReducedInhomBayesianSys} now has the reduced impulse response
\begin{align}\label{eq:ReducedImpulseResponse}
    \mbf{h}_r(t)= \mbf{\Gamma}_\epsilon^{-1/2}\mbf{C}_re^{\mbf{A}_rt}\mbf{L}_{\pr,r}.
\end{align}
For the remainder of this paper, we assume that the full prior-driven system \eqref{eq:PriorDrivenSys} is stable and \textit{minimal} in the sense that it has only non-zero Hankel singular values, in which case $\mbf{\Sigma}$ is invertible. Given this minimality assumption, the bases $\mbf{W}:=\mbf{W}_d = \mbf{L}\mbf{Z}\mbf{\Sigma}^{-1/2}\in\R^{d\times d}$ and $\mbf{V}:=\mbf{V}_d =\mbf{R}\mbf{U}\mbf{\Sigma}^{-1/2}\in\R^{d\times d}$ are well-defined. In particular, the PD-BT-reduced bases $\mathbf{W}_r,\mathbf{V}_r\in\R^{d\times r}$ that we defined above are given by the first $r$ columns of $\mbf{W}$ and $\mbf{V}$, respectively, for a truncation rank $r\leq d$.

Note that, even if $\PrCov$ is rank-deficient, the system \eqref{eq:PriorDrivenSys} can be minimal. The minimality assumption is not very restrictive and is commonly used in systems theory, see, e.g., \cite{Antoulas2005Book,BEATTIE2017inhom,Moore1981BT,Qian2021Balancing}. On the other hand, for small Hankel singular values, computing the full bases $\mbf{V}$ and $\mbf{W}$ is numerically ill-conditioned and is usually avoided. For constructing the reduced system \eqref{eq:ReducedInhomBayesianSys} numerically, only the reduced bases $\mbf{V}_r$ and $\mbf{W}_r$ are required.

Given the minimality assumption, for a reduced system constructed by PD-BT, the difference between the full impulse response $\mbf{h}(t)$ \eqref{eq:FullImpulseResponse} and the reduced impulse response $\mbf{h}_r(t)$ \eqref{eq:ReducedImpulseResponse} can be bounded as follows:
\begin{proposition}[{\cite[Theorem 3.1 and 3.2]{BEATTIE2017inhom}}]
\label{prop:SplittingErrorBound}
    Consider the full system~\eqref{eq:PriorDrivenSys} and its rank-$r$ reduced version~\eqref{eq:ReducedInhomBayesianSys} obtained by PD-BT. Let $\mbf{S}$ be a solution of the Sylvester equation $$\mbf{W}^\top\mbf{A}^\top\mbf{V}\mbf{S} + \mbf{SA}_r + \mbf{W}^\top\mbf{C}^\top\mbf{\Gamma}_\epsilon^{-1}\mbf{C}_r =\mbf{0},$$ and let $\bar{\mbf{S}}\in\R^{(d-r)\times r}$ collect the last $d-r$ rows of $\mbf{S}$. Similarly, let $\bar{\mbf{V}},\bar{\mbf{W}}\in\R^{d\times (d-r)}$ collect the last $d-r$ columns of $\mbf{V},\mbf{W}$ and let $\bar{\mbf{\Sigma}}\coloneqq\mathrm{diag}(\sigma_{r+1},\ldots,\sigma_d)\in\R^{(d-r)\times (d-r)}$ be the diagonal matrix of the $d-r$ truncated HSVs of $\mathbf{P}\mbf{Q}$. Let $\bar{\mbf L}_\pr\coloneqq\bar{\mbf{V}}^\top\mbf{L}_\pr\in\R^{(d-r)\times s}$ and  $\bar{\mbf{A}}\coloneqq\mbf{V}_r^\top\mbf{A}\bar{\mbf{W}}\in\R^{r\times (d-r)}$.
        Then, for the impulse responses $\mbf{h}(t)$ of the full system~\eqref{eq:PriorDrivenSys} and $\mbf{h}_r(t)$ of the balanced, rank-$r$ reduced system~\eqref{eq:ReducedInhomBayesianSys},
       \begin{align}\label{eq:InhomImpulseErrorBound}
            \| \mbf{h}-\mbf{h}_r\|^2_{L_2} \leq \mathrm{trace}\big[ \big(\bar{\mbf L}_\pr\bar{\mbf L}_\pr^\top+2\bar{\mbf{S}}\bar{\mbf{A}}\big)\bar{\mbf{\Sigma}}\big].
        \end{align}
\end{proposition}
\Cref{prop:SplittingErrorBound} highlights that the quality of the approximation of the full impulse response by its PD-BT-reduced version depends on the matrix $\bar{\mbf{\Sigma}}$  of truncated HSVs of the prior-driven system~\eqref{eq:PriorDrivenSys}.

For a reduced system constructed by PD-TLBT, the difference between the full impulse response $\mbf{h}(t)$ \eqref{eq:FullImpulseResponse} and the reduced impulse response $\mbf{h}_r(t)$ \eqref{eq:ReducedImpulseResponse} can be bounded as follows:
\begin{proposition}\label[proposition]{prop:TLSplittingErrorBound}
    Consider the full system~\eqref{eq:PriorDrivenSys} and its rank-$r$ reduced version~\eqref{eq:ReducedInhomBayesianSys} obtained by PD-TLBT. Let $\mathbf{P}^{(\textup{red})}_{\mathcal{T},r} \coloneqq\int_0^\mathcal{T} e^{\mathbf{A}_r t}\mbf{L}_{\pr,r}\mbf{L}_{\pr,r}^\top e^{\mathbf{A}_r^\top t}\,\text{d}t$, and let $\mathbf{P}^{(\textup{mix})}_{\mathcal{T},r} \coloneqq\int_0^\mathcal{T} e^{\mathbf{A} t}\mbf{L}_\pr\mbf{L}_{\pr,r}^\top e^{\mathbf{A}_r^\top t}\,\text{d}t$.
    Then, for the impulse responses $\mbf{h}(t)$ of the full system~\eqref{eq:PriorDrivenSys} and $\mbf{h}_r(t)$ of the balanced, rank-$r$ reduced system~\eqref{eq:ReducedInhomBayesianSys},
       \begin{align}\label{eq:TLInhomImpulseErrorBound}
            \| &\mbf{h}-\mbf{h}_r\|^2_{L_{2,\mathcal{T}}} \\
            &\leq \mathrm{trace}\big[\mbf{\Gamma}_\epsilon^{-1}\mbf{C}\mathbf{P}_{\mathcal{T}}\mbf{C}^\top\big]+\mathrm{trace}\big[\mbf{\Gamma}_\epsilon^{-1}\mbf{C}_r\mathbf{P}^{(\textup{red})}_{\mathcal{T},r}\mbf{C}_r^\top\big]-2\,\mathrm{trace}\big[\mbf{\Gamma}_\epsilon^{-1}\mbf{C}\mathbf{P}^{(\textup{mix})}_{\mathcal{T},r}\mbf{C}_r^\top\big].\notag      
        \end{align}
        \begin{proof}
        Using the definitions of the time-limited $L_2$-norm \eqref{eq:TLL2Norm} and the full and approximate impulse responses \eqref{eq:FullImpulseResponse} and \eqref{eq:ReducedImpulseResponse} we proceed as in \cite{Redmann2018H2TLBT} and obtain:
        \begin{align*}
            \| \mbf{h}&-\mbf{h}_r\|^2_{L_{2,\mathcal{T}}}\\
            \overset{\eqref{eq:TLL2Norm}}{=}&\int_0^\mathcal{T}\norm{\mbf{h}(t)-\mbf{h}_r(t)}_F^2\,\text{d}t\\
            \overset{\eqref{eq:FullImpulseResponse},\, \eqref{eq:ReducedImpulseResponse}}{=} &\int_0^\mathcal{T}\norm{ \mbf{\Gamma}_\epsilon^{-1/2}\mbf{C}e^{\mbf{A}t}\mbf{L}_\pr- \mbf{\Gamma}_\epsilon^{-1/2}\mbf{C}_re^{\mbf{A}_rt}\mbf{L}_{\pr,r}}_F^2\,\text{d}t\\
            \overset{\text{Def. } \Norm{\cdot}_F}{=}&\int_0^\mathcal{T}\mathrm{trace}\big[ \mbf{\Gamma}_\epsilon^{-1/2}\mbf{C}e^{\mathbf{A} t}\PrCov e^{\mathbf{A}^\top t}\mbf{C}^\top \mbf{\Gamma}_\epsilon^{-1/2}\big]\,\text{d}t\\
            +&\int_0^\mathcal{T}\mathrm{trace}\big[ \mbf{\Gamma}_\epsilon^{-1/2}\mbf{C}_re^{\mathbf{A}_r t}\mbf{L}_{\pr,r}\mbf{L}_{\pr,r}^\top e^{\mathbf{A}_r^\top t}\mbf{C}_r^\top \mbf{\Gamma}_\epsilon^{-1/2}\big]\,\text{d}t\\
            -2\,&\int_0^\mathcal{T}\mathrm{trace}\big[ \mbf{\Gamma}_\epsilon^{-1/2}\mbf{C}e^{\mathbf{A} t}\mbf{L}_\pr\mbf{L}_{\pr,r}^\top e^{\mathbf{A}_r^\top t}\mbf{C}_r^\top \mbf{\Gamma}_\epsilon^{-1/2}\big]\,\text{d}t\\
            =&\enspace\mathrm{trace}\big[\mbf{\Gamma}_\epsilon^{-1}\mbf{C}\mathbf{P}_{\mathcal{T}}\mbf{C}^\top\big]+\mathrm{trace}\big[\mbf{\Gamma}_\epsilon^{-1}\mbf{C}_r\mathbf{P}^{(\textup{red})}_{\mathcal{T},r}\mbf{C}_r^\top\big]-2\,\mathrm{trace}\big[\mbf{\Gamma}_\epsilon^{-1}\mbf{C}\mathbf{P}^{(\textup{mix})}_{\mathcal{T},r}\mbf{C}_r^\top\big],
        \end{align*}
        where the last equality results from the linear and cyclic properties of the trace and from inserting the definitions of $\mathbf{P}_{\mathcal{T}}$ from \eqref{eq:TL_Gramians}, $\mathbf{P}^{(\textup{red})}_{\mathcal{T},r}$ and $\mathbf{P}^{(\textup{mix})}_{\mathcal{T},r}$.
        \end{proof}
\end{proposition}
\begin{remark}
\label{remark:TLBT_bound_HSVs}
    The bound \eqref{eq:TLInhomImpulseErrorBound} can also be expressed in terms of the truncated HSVs $\bar{\mbf{\Sigma}}$ of the prior-driven system~\eqref{eq:PriorDrivenSys}, but in a more complicated and computationally impractical way. For more details, see \cite[Theorem 2.3]{Redmann2018H2TLBT}.
\end{remark}
Note that $\norm{\mbf{h}(t_k)-\mbf{h}_r(t_k)}_F^2$ is non-negative and finite for all times $t_k$. Hence, for each $k=1,\ldots,n$ there is a constant $C_k\geq 0$ depending on $\norm{\mbf{h}(t)-\mbf{h}_r(t)}_F$ and on the observation times $t_k$ such that $\norm{\mbf{h}(t_k)-\mbf{h}_r(t_k)}_F^2 \leq C_k \int_{t_{k-1}}^{t_k} \norm{\mbf{h}(t)-\mbf{h}_r(t)}_F^2 \mathrm{d}t$, where we set $t_0 = 0$. We set $\kappa \coloneqq \max_k C_k$. Then, by the proof of \cite[Theorem 3.6]{2025PD-BT}, for observation times $t_n\leq\mathcal{T}$ in \eqref{eq:BayesianSystem_noise}, the following inequality holds:
\begin{align}\label{eq:integral_limit}
 \begin{split}
    \sum_{k=1}^n\|\mbf{h}(t_k)-\mbf{h}_r(t_k)\|_F^2 \leq &\kappa \int_{0}^{t_n} \|\mbf{h}(t)-\mbf{h}_r(t)\|_F^2 \mathrm{d}t \\
    \overset{t_n\leq\mathcal{T} \text{ and }\eqref{eq:TLL2Norm}}{\leq}&\kappa\| \mbf{h}-\mbf{h}_r\|^2_{L_{2,\mathcal{T}}}\overset{\eqref{eq:L2Norm}}{\leq} \kappa\| \mbf{h}-\mbf{h}_r\|^2_{L_2}.     
 \end{split}
\end{align}
\begin{remark}
The constant $\kappa$ as defined above is inversely proportional to the time step between measurements and thus becomes arbitrarily large in the limit of continuous-time measurements.
      However, for this work it suffices that $\kappa$ is finite and constant with respect to $r$, which is the case provided that there are only finitely many observations in time. These properties of $\kappa$ suffice because the purpose of the error bound \eqref{eq:integral_limit} is to identify trends in the error associated to the reduced order model, and not to determine the error with high accuracy \cite[Remark 3.7]{2025PD-BT}. 
      In particular, \eqref{eq:integral_limit} and \Cref{prop:SplittingErrorBound,prop:TLSplittingErrorBound} allow one to relate the decay in the sum of impulse response errors, evaluated at the times of measurement, to the decay of the truncated HSVs.
      This is possible even for the case of finitely many high-frequency observations, in which case $\kappa$ is large but remains finite.
\end{remark}
The reduced prior-driven system \eqref{eq:ReducedInhomBayesianSys} now serves for MOR in Bayesian inference by using the reduced system matrices $\mbf{A}_r,\mbf{C}_r$ and the reduced left basis $\mbf{V}_r$ to define a reduced forward map
\begin{align}\label{eq:PD-BT forward map}
    \widehat{\fwdmap} = \begin{bmatrix}
        \mbf{C}_re^{\mbf{A}_rt_1}\\ \cdots\\ \mbf{C}_re^{\mbf{A}_rt_n}
    \end{bmatrix}\mbf{V}_r^\top\in \R^{d_{\rm obs}\times d},
\end{align}
which has the same dimensions as the full forward map $\mbf{G}$ \eqref{eq:BLIPproperties} but requires matrix products of much smaller dimension $r\ll d$. We can then apply $\widehat{\fwdmap}$ in \eqref{eq:PD-BT PosBoth} to obtain an approximate posterior mean and covariance.

Up to now, only a bound between the noise-free outputs  of~\eqref{eq:BayesianSystem_nonoise}, namely \linebreak\mbox{$\mathbf{Y} \coloneqq \begin{bmatrix}\mathbf{y}(t_1)^\top,\ldots,\mathbf{y}(t_n)^\top\end{bmatrix}^\top \in \mathbb{R}^{d_\mathrm{obs}}$}, and the noise-free reduced outputs \linebreak\mbox{$\mathbf{Y}_r \coloneqq \begin{bmatrix}\mathbf{y}_r(t_1)^\top,\ldots,\mathbf{y}_r(t_n)^\top\end{bmatrix}^\top \in \mathbb{R}^{d_\mathrm{obs}}$}, where $\mathbf{y}_r(t)\coloneqq\mathbf{C}_re^{\mathbf{A}_rt}\mathbf{V}_r^\top\mathbf{p}$ for every $t$, was given in \cite[Theorem 3.6]{2025PD-BT}. In \cref{sec:PD_BT_bound} we shall present bounds on $\|\PosCov - \widehat{\mbf{\Gamma}}_\pos\|_F$ and $\|\PosMn(\data)- \widehat{\pmb{\mu}}_\pos(\data)\|_2$ with $\widehat{\mbf{\Gamma}}_\pos$ and $\widehat{\pmb{\mu}}_\pos(\data)$ from \eqref{eq:PD-BT PosBoth} obtained by PD-(TL)BT. Unlike the bound on the output approximation, the bounds that we will present precisely describe the quality of the approximation in the posterior, which is the solution to the Gaussian inverse problem. We will demonstrate that the quality of the posterior approximation depends on the matrix $\bar{\mbf{\Sigma}}$ of truncated HSVs of the prior-driven system~\eqref{eq:PriorDrivenSys}. To derive this bound, we will rely on local Lipschitz stability of Gaussian posteriors, which we introduce next.

\section{Local Lipschitz stability of Gaussian posteriors}
\label{sec:Lipschitz bounds}

In this section, we show in \Cref{theorem03} first that the posterior covariance is a locally Lipschitz continuous function of the forward model, and then use this to show the same is true of the posterior mean.
Local Lipschitz stability results of this type are proven for possibly non-linear and non-Gaussian Bayesian inverse problems in \cite{Sprungk2020LocalLipschitz}, where the error in the posterior measure is controlled.
Our focus on linear Gaussian inverse problems implies that it suffices to control the error in the posterior mean and covariance.
For linear Gaussian inverse problems, Theorem 3.3 and Theorem 3.4 of \cite{SanzAlonso2024} provide bounds on the error of the approximate posterior mean and covariance, respectively, where the approximations are constructed using projections, and the approximation errors are controlled by operator norm bounds of the projection errors \cite[Assumption 3.1]{SanzAlonso2024}.
For \Cref{theorem03} below, we do not assume that the approximations are obtained by projections; we only assume that $\fwdmap$ and $\hatFwdmap$ are linear.

For \Cref{theorem03} below, recall that for a matrix $\mbf{A}$, $\norm{\mbf{A}}_p$ denotes the $p$-Schatten norm of $\mbf{A}$; see \Cref{section_notation}.

\begin{restatable}{theorem}{errorBoundppHessianPerturbation}
  \label{theorem03}
  Let $\data$, $\PrCov$, $\PrMn$, $\ObsCov$ and $\fwdmap$ satisfy \Cref{assumption_main}, and let $1\leq p\leq \infty$. 
Let $\hatFwdmap$ be another linear forward model and let $\mbf{L}_\pr$ be a possibly non-symmetric square root of $\PrCov$, i.e. $\mbf{L}_\pr\mbf{L}_\pr^\top=\PrCov$.
Then $\PosCov$ and $\hatPosCov$ in \eqref{eq:pos_covariance} and \eqref{eq:pos_covariance_approx} satisfy, for some $0<C=C(\fwdmap,\hatFwdmap,\PrCov,\ObsCov,\mbf{L}_\pr)<\infty$,
   \begin{align}
 \label{eq:error_covariance_for_perturbed_forward_model_by_ppH}
 \Norm{ \PosCov -\hatPosCov}_p\leq  C\Norm{ \ObsCov^{-1/2}(\fwdmap-\hatFwdmap )\mbf{L}_\pr}_p.
 \end{align}
If in addition $\PrMn\in\range \mbf{L}_{\pr}$, then $\PosMn(\data)$ and $\hatPosMn(\data)$ in \eqref{eq:pos_mean} and \eqref{eq:pos_mean_approx} satisfy, for some $0<C'=C'(C,\fwdmap,\hatFwdmap,\PrCov,\ObsCov,\mbf{L}_\pr,\PrMn,\data)<\infty$,
\begin{align}
 \label{eq:error_mean_for_perturbed_forward_model_by_ppH}
 \Norm{\PosMn(\data)-\hatPosMn(\data)}_2\leq  C' \Norm{ \ObsCov ^{-1/2} (\fwdmap-\hatFwdmap) \mbf{L}_\pr}_\infty.
\end{align}
 \end{restatable}
We prove \Cref{theorem03} in \Cref{section_proofs_local_lipschitz_bounds}, and define the scalars $C$ and $C'$ in \eqref{eq:error_covariance_for_perturbed_forward_model_by_ppH_constant} and \eqref{eq:error_mean_for_perturbed_forward_model_by_ppH_constant}, respectively.

 The significance of \Cref{theorem03} is to show that, for finite-dimensional linear Gaussian inverse problems, the posterior covariance and mean are locally Lipschitz continuous functions of the square roots of the prior-preconditioned Hessian corresponding to the linear forward model.
 Alternatively, the errors in the approximate posterior covariance and mean are of first order in the errors in the square root of the prior-preconditioned Hessian.
 
 Since the prior-preconditioned Hessian includes the prior covariance and observation precision, the bounds incorporate regularization by both the prior and the noise: larger (respectively, smaller) prior covariances and larger (resp. smaller) observation precisions lead to larger (resp. smaller) errors.
 In the numerical results that we report in \cref{subsec:Numerics discussion}, we shall show that this dependence on the prior is not only a feature of the error bounds, but also of the actual errors.

The hypothesis that $\PrMn\in\range \mbf{L}_\pr$ in \Cref{theorem03} is not new; see \cite[Example 6.23]{Stuart2010}.
The hypothesis holds in the setting described in \cref{subsec:LTI FOM}, where we assume $\PrMn=\mbf{0}$.
If $\PrMn\neq \mbf{0}$ and $\PrMn\in\range \PrCov^{1/2}$, then for every $\mbf{L}_\pr$ such that $\mbf{L}_\pr\mbf{L}_\pr^\top=\PrCov$, $\PrMn\in\range \mbf{L}_\pr$ holds.
This is because by the singular value decomposition and the finite-dimensional setting, $\range \PrCov^{1/2}=\range \PrCov$, and because $\PrCov=\mbf{L}_\pr\mbf{L}_\pr^\top$ implies $\range \PrCov\subseteq \range \mbf{L}_\pr$.
In particular, if $\PrCov\succ 0$, then $\range \PrCov^{1/2}=\R^d$, and the hypothesis holds.
We conjecture that the condition $\PrMn\in\range\mbf{L}_\pr$ in \Cref{theorem03} is not only sufficient but also necessary for an error bound of the form \eqref{eq:error_mean_for_perturbed_forward_model_by_ppH}. 

We conclude this section with a discussion of the case of nonlinear inverse problems. If the forward map $\fwdmap$ is nonlinear, then the posterior is in general not Gaussian and \eqref{eq:Posterior} does not hold, even for Gaussian priors and observation noise. In this setting, one can construct a Gaussian approximation of the non-Gaussian reference posterior by using a linear approximation $\hatFwdmap$ of $\fwdmap$ such as the `Gauss--Newton Hessian' or the Laplace approximation; see e.g. \cite{BuiThanh2012} and \cite{Kretschmann2022} respectively. 
In the most general setting of possibly nonlinear $\fwdmap$ and non-Gaussian priors and observation noise, the Bayesian approach to inverse problems remains valid under mild integrability and continuity conditions \cite[Section 4.4]{Stuart2010}. 
Moreover, the error in the approximate non-Gaussian posterior in various metrics and information-theoretic divergences can be bounded in terms of the error $\hatFwdmap-\fwdmap$ for nonlinear forward models and non-Gaussian priors using the local Lipschitz stability of Bayesian inverse problems; see e.g. \cite[Remark 4]{Sprungk2020LocalLipschitz}, and see \cite{CvetkovicLie2025} for a more recent analysis of local Lipschitz stability.
As with \Cref{theorem03}, the purpose of the bounds in \cite{CvetkovicLie2025,Sprungk2020LocalLipschitz,Stuart2010} is not to provide tight control of the error, but rather to indicate how the error of the approximate posterior depends on errors such as $\fwdmap-\hatFwdmap$.

\section{Posterior error bounds for prior-driven balancing}\label{sec:PD_BT_bound}
Recall that the posterior approximation \eqref{eq:PD-BT PosBoth} obtained using PD-(TL)BT is defined by the reduced forward map $\hatFwdmap$ from \eqref{eq:PD-BT forward map} associated to the prior-driven homogeneous LTI system \eqref{eq:PriorDrivenSys}, whereas the exact posterior \eqref{eq:Posterior} is defined by the exact forward map $\fwdmap$ from \eqref{eq:BLIPproperties} associated to the LTI system \eqref{eq:BayesianSystem}. We will now use the local Lipschitz bound for a perturbed square root arising from a perturbed forward model in the particular case $\PrMn=\mbf{0}$, given by~\eqref{eq:error_covariance_for_perturbed_forward_model_by_ppH} and~\eqref{eq:error_mean_for_perturbed_forward_model_by_ppH}, to derive error bounds for the posterior approximation.

It is crucial to note that, for our described setting, the rows of the square root of the prior-preconditioned Hessian are exactly the values of the impulse response function \eqref{eq:FullImpulseResponse} at the observation times $t_1,\ldots,t_n$. Similarly, for the reduced system, the rows of the square root of the corresponding reduced prior-preconditioned Hessian are the values of reduced impulse response function \eqref{eq:ReducedImpulseResponse} at $t_1,\ldots,t_n$. By the definitions of $\ObsCov$ and $\mbf{L}_{\pr,r}$ in \eqref{eq:BLIPproperties} and below \eqref{eq:ReducedInhomBayesianSys} respectively, we obtain
\begin{align}\label{eq:PD-BT_reasoning}
\begin{split}
        &\ObsCov^{-1/2}(\mbf{G-\widehat{G}}) \mbf{L}_\pr\\
        &= \begin{bmatrix} \mbf{\Gamma}_\epsilon^{-1/2}\\ & \ddots \\ && \mbf{\Gamma}_\epsilon^{-1/2} \end{bmatrix}\cdot\Bigg(\begin{bmatrix}\mbf{C} e^{\mbf{A}t_1}\\ \vdots\\\mbf{C} e^{\mbf{A}t_n}\end{bmatrix}-\begin{bmatrix}\mbf{C}_r e^{\mbf{A}_rt_1}\\ \vdots\\\mbf{C}_r e^{\mbf{A}_rt_n}\end{bmatrix}\mbf{V}_r^\top\Bigg) \cdot \begin{bmatrix}\mbf{L}_\pr \end{bmatrix}\\
        &= \begin{bmatrix} \mbf{\Gamma}_\epsilon^{-1/2}\mbf{C} e^{\mbf{A}t_1}\mbf{L}_\pr\\ \vdots\\  \mbf{\Gamma}_\epsilon^{-1/2}\mbf{C} e^{\mbf{A}t_n}\mbf{L}_\pr\end{bmatrix}-\begin{bmatrix} \mbf{\Gamma}_\epsilon^{-1/2}\mbf{C}_r e^{\mbf{A}_rt_1}\mbf{L}_{\pr,r}\\ \vdots\\  \mbf{\Gamma}_\epsilon^{-1/2}\mbf{C}_r e^{\mbf{A}_rt_n}\mbf{L}_{\pr,r}\end{bmatrix} = \begin{bmatrix}\mbf{h}(t_1)-\mbf{h}_r(t_1)\\\vdots\\\mbf{h}(t_n)-\mbf{h}_r(t_n) \end{bmatrix}.    
\end{split}
\end{align}
The impulse response of a system is a key system-theoretic property that completely defines the behavior of an LTI system \cite[Chapter 4.1]{Antoulas2005Book}. Meanwhile, the prior-preconditioned Hessian is crucial for analyzing linear Gaussian inverse problems \eqref{eq:generic_IVP}. For example, the spectral decay of the prior-preconditioned Hessian controls the quality of certain optimal low-rank approximations of the posterior mean and covariance \cite{CarereLie2025a,CarereLie2025b,Flath2011,Spantini2015Optimal}. Thus, \eqref{eq:PD-BT_reasoning} reveals a novel connection between Bayesian inference and systems theory. This motivates the use of system-theoretic MOR of the prior-driven system \eqref{eq:PriorDrivenSys} to reduce the computational burden of solving Bayesian smoothing problems for LTI systems.

The above-mentioned connection between Bayesian inference and systems theory allows us to bound the error in the approximation of the square root of the prior-preconditioned Hessian by the truncated HSVs of the prior-driven system:
    \begin{align}\label{eq:PriPreHessBound}
       \|\ObsCov^{-1/2}(\mbf{G-\widehat{G}}) \mbf{L}_\pr\|_F^2&\overset{\eqref{eq:PD-BT_reasoning}}{=}\left\lVert\begin{bmatrix}\mbf{h}(t_1)-\mbf{h}_r(t_1)\\\vdots\\\mbf{h}(t_n)-\mbf{h}_r(t_n) \end{bmatrix}\right\rVert^2_F\notag \\
       &= \sum_{i=1}^n \|\mbf{h}(t_i)-\mbf{h}_r(t_i)\|^2_F \notag\\
       &\overset{\eqref{eq:InhomImpulseErrorBound},\eqref{eq:integral_limit}}{\leq}\kappa \cdot \mathrm{trace}\big[\big(\bar{\mbf L}_\pr\bar{\mbf L}_\pr^\top+2\bar{\mbf{S}}\bar{\mbf{A}}\big)\bar{\mbf{\Sigma}} \big].
    \end{align}
The above observations then lead to the following bounds on the posterior mean and covariance:
\begin{theorem}[PD-BT posterior bound]\label{thm:PD-BT bounds}
Let $\widehat{\pmb{\mu}}_\pos(\data)$ be the approximate mean and $\widehat{\mbf{\Gamma}}_\pos$ the approximate covariance~\eqref{eq:PD-BT PosBoth} of the Bayesian smoothing problem~\eqref{eq:BayesianSystem} obtained by PD-BT. Let $\bar{\mbf L}_\pr$, $\bar{\mbf{S}}$, $\bar{\mbf{A}}$ and $\bar{\mbf{\Sigma}}$ be as in \cref{prop:SplittingErrorBound} and $\kappa$ be as in~\eqref{eq:integral_limit} and the scalars $C$ and $C'$ be as in \eqref{eq:error_covariance_for_perturbed_forward_model_by_ppH} and \eqref{eq:error_mean_for_perturbed_forward_model_by_ppH}, respectively. Then, the difference between the full and the reduced posterior quantities can be bounded as follows:
\begin{align}\label{eq:Frobenius PosCov bound}
\begin{split}
\|\PosCov - \widehat{\mbf{\Gamma}}_\pos\|_F \leq  C \sqrt{\kappa \cdot \mathrm{trace}\big[\big(\bar{\mbf L}_\pr\bar{\mbf L}_\pr^\top+2\bar{\mbf{S}}\bar{\mbf{A}}\big)\bar{\mbf{\Sigma}} \big]},
\end{split}
\end{align}
and
\begin{align}\label{eq:Frobenius PosMean bound}
    \|\PosMn(\data)- \widehat{\pmb{\mu}}_\pos(\data)\|_2 = \|&\PosMn(\data)- \widehat{\pmb{\mu}}_\pos(\data)\|_F\\
    &\leq C'\cdot\sqrt{\kappa \cdot \mathrm{trace}\big[\big(\bar{\mbf L}_\pr\bar{\mbf L}_\pr^\top+2\bar{\mbf{S}}\bar{\mbf{A}}\big)\bar{\mbf{\Sigma}} \big]}.\notag
\end{align}
\end{theorem}
\begin{proof}
    The argument follows by applying the local Lipschitz bounds \eqref{eq:error_covariance_for_perturbed_forward_model_by_ppH} for the covariance and \eqref{eq:error_mean_for_perturbed_forward_model_by_ppH} for the mean from \Cref{theorem03} with $p=2$, together with the bound~\eqref{eq:PriPreHessBound} between the full and the reduced prior-preconditioned Hessian.
\end{proof}

Similar considerations can be made for PD-TLBT:
\begin{theorem}[PD-TLBT posterior bound]\label{thm:PD-TLBT bounds}
Consider the same setting as in \cref{thm:PD-BT bounds} with the approximate posterior mean $\widehat{\pmb{\mu}}_\pos(\data)$ and covariance $\widehat{\mbf{\Gamma}}_\pos$ obtained by PD-TLBT. Let $\mathbf{P}^{(\textup{red})}_{\mathcal{T},r}$, $\mathbf{P}^{(\textup{mix})}_{\mathcal{T},r}$ be as in \cref{prop:TLSplittingErrorBound} and $\kappa$ be as in~\eqref{eq:integral_limit} and the scalars $C$ and $C'$ be as in \eqref{eq:error_covariance_for_perturbed_forward_model_by_ppH} and \eqref{eq:error_mean_for_perturbed_forward_model_by_ppH}, respectively. Then, the difference between the full and the reduced posterior quantities can be bounded as follows:
\begin{align}\label{eq:TL Frobenius PosCov bound}
\begin{split}
\|&\PosCov - \widehat{\mbf{\Gamma}}_\pos\|_F \leq  C \cdot \\&\sqrt{\kappa \cdot\Big( \mathrm{trace}\big[ \mbf{\Gamma}_\epsilon^{-1}\mbf{C}\mathbf{P}_{\mathcal{T}}\mbf{C}^\top\big]+\mathrm{trace}\big[ \mbf{\Gamma}_\epsilon^{-1}\mbf{C}_r\mathbf{P}^{(\textup{red})}_{\mathcal{T},r}\mbf{C}_r^\top\big]-2\,\mathrm{trace}\big[ \mbf{\Gamma}_\epsilon^{-1}\mbf{C}\mathbf{P}^{(\textup{mix})}_{\mathcal{T},r}\mbf{C}_r^\top\big]\Big)},
\end{split}
\end{align}
and
\begin{align}\label{eq:TL Frobenius PosMean bound}
\begin{split}
    \|&\PosMn(\data)- \widehat{\pmb{\mu}}_\pos(\data)\|_2 = \|\PosMn(\data)- \widehat{\pmb{\mu}}_\pos(\data)\|_F\leq C'\cdot 
    \\&\sqrt{\kappa \cdot \Big(\mathrm{trace}\big[ \mbf{\Gamma}_\epsilon^{-1}\mbf{C}\mathbf{P}_{\mathcal{T}}\mbf{C}^\top\big]+\mathrm{trace}\big[ \mbf{\Gamma}_\epsilon^{-1}\mbf{C}_r\mathbf{P}^{(\textup{red})}_{\mathcal{T},r}\mbf{C}_r^\top\big]-2\,\mathrm{trace}\big[ \mbf{\Gamma}_\epsilon^{-1}\mbf{C}\mathbf{P}^{(\textup{mix})}_{\mathcal{T},r}\mbf{C}_r^\top\big]\Big)}.
\end{split}
\end{align}
\end{theorem}
\begin{proof}
    The argument follows by applying the local Lipschitz bounds \eqref{eq:error_covariance_for_perturbed_forward_model_by_ppH} for the covariance and \eqref{eq:error_mean_for_perturbed_forward_model_by_ppH} for the mean from \Cref{theorem03} with $p=2$, together with the following bound between the full and the PD-TLBT-reduced prior-preconditioned Hessian:
    \begin{align*}
       \|\ObsCov^{-1/2}(\mbf{G-\widehat{G}}) \mbf{L}_\pr\|_F^2
       \overset{\eqref{eq:PD-BT_reasoning}}{=} &\sum_{i=1}^n \|\mbf{h}(t_i)-\mbf{h}_r(t_i)\|^2_F \notag\\
       \overset{\eqref{eq:TLInhomImpulseErrorBound},\eqref{eq:integral_limit}}{\leq}&\kappa \cdot \mathrm{trace}\big[\mbf{\Gamma}_\epsilon^{-1}\mbf{C}\mathbf{P}_{\mathcal{T}}\mbf{C}^\top\big]+\mathrm{trace}\big[\mbf{\Gamma}_\epsilon^{-1}\mbf{C}_r\mathbf{P}^{(\textup{red})}_{\mathcal{T},r}\mbf{C}_r^\top\big]\\&-2\,\mathrm{trace}\big[\mbf{\Gamma}_\epsilon^{-1}\mbf{C}\mathbf{P}^{(\textup{mix})}_{\mathcal{T},r}\mbf{C}_r^\top\big].
    \end{align*}
\end{proof}
\Cref{thm:PD-BT bounds,thm:PD-TLBT bounds} imply that the posterior approximation error in the Bayesian smoothing problem \cref{eq:BayesianSystem} depends on the truncated HSVs $\bar{\mbf{\Sigma}}$ of the prior-driven system \eqref{eq:PriorDrivenSys}, and thus on a system-theoretic invariant; see \Cref{remark:TLBT_bound_HSVs}. In particular, small truncated HSVs lead to small bounds \eqref{eq:Frobenius PosCov bound}--\eqref{eq:TL Frobenius PosMean bound}, and thus guarantee good posterior approximation. Note that the HSVs only enter through their square root, so that the convergence of $\|\PosCov - \widehat{\mbf{\Gamma}}_\pos\|_F$ and $\|\PosMn(\data)- \widehat{\pmb{\mu}}_\pos(\data)\|_2$ to zero is not very fast.

Recall in the discussion after \Cref{theorem03} that the error bounds \eqref{eq:error_covariance_for_perturbed_forward_model_by_ppH} and \eqref{eq:error_mean_for_perturbed_forward_model_by_ppH} on the approximate posterior covariance and mean incorporate regularization by the prior and the noise.
In \eqref{eq:Frobenius PosCov bound} and \eqref{eq:Frobenius PosMean bound}, the dependence on $\PrCov$ appears via the projected square root prior covariance $\bar{\mbf L}_\pr$, while the dependence on $\ObsCov^{-1}$ appears via the row submatrix $\bar{\mbf{S}}$ of the solution $\mbf{S}$ to a Sylvester equation that depends on $ \mbf{\Gamma}_\epsilon^{-1}$; see \Cref{prop:SplittingErrorBound}, and recall the definition of $\ObsCov$ in \eqref{eq:BLIPproperties}.
In \eqref{eq:TL Frobenius PosCov bound} and \eqref{eq:TL Frobenius PosMean bound}, the dependence on $\PrCov$ appears indirectly via the Gramian $\mathbf{P}_\mathcal{T}$ in \eqref{eq:TL_Gramians} and its rank-$r$ variants $\mathbf{P}^{(\textup{red})}_{\mathcal{T},r}$ and $\mathbf{P}^{(\textup{mix})}_{\mathcal{T},r}$ defined in \Cref{prop:TLSplittingErrorBound}, while the dependence on $\ObsCov$ appears explicitly via the $ \mbf{\Gamma}_\epsilon$ terms.
In the following section, we will present results that illustrate the dependence of the errors and their bounds on the prior.

\section{Numerical experiments}\label{sec:Numerics}
We now present numerical experiments comparing the approximation error in the posterior mean and covariance for Bayesian inference of the initial condition of an LTI system with the respective bounds derived in \cref{sec:PD_BT_bound}. We describe the benchmark test problems considered in \cref{subsec:test problem} and discuss the experimental results in \cref{subsec:Numerics discussion}.

\subsection{Test problems}\label{subsec:test problem}
\paragraph{ISS1R benchmark}
We consider the same benchmark example as in \cite[Sec. 4]{2025PD-BT}: The ISS1R benchmark model represents the structural dynamics of the Russian service module of the International Space Station. It has three inputs: the roll, pitch, and yaw jets, and three outputs: the respective roll, pitch, and yaw gyroscope measurements \cite{Antoulas2005Book}. Through discretization, a stable LTI system with system matrices of the following dimensions is obtained: $\mathbf{A}\in \mathbb{R}^{270 \times 270}$, $\mathbf{B}\in \mathbb{R}^{270 \times 3}$, and $\mathbf{C}\in \mathbb{R}^{3 \times 270}$. The system matrices for the benchmark are available at \url{http://slicot.org/20-site/126-benchmark-examples-for-model-reduction}.

To generate the vector $\data$ of noisy output measurements, we initialized the system at a reference initial condition $\mathbf{p}$ drawn from the prior distribution $\mathcal{N}(\mathbf{0},\PrCov)$, and then simulated the linear system from time $t_0 = 0$ up to the final time of measurement $t_n = 8$. Measurements were computed according to \eqref{eq:BayesianSystem_noise} using the resulting trajectory, equidistant times $t_k=k\Delta t $ for $\Delta t=0.1$, and Gaussian measurement noise $\mathcal{N}(\mathbf{0},\mathbf{\Gamma}_\epsilon)$ for $\mathbf{\Gamma}_\epsilon=\mathrm{diag}(0.0025^2,0.0005^2,0.0005^2)$. This is the same choice of $\mathbf{\Gamma}_\epsilon$ as in \cite{Koe2022TLBTDA,2025PD-BT,Qian2021Balancing}.

In the discussion after the proof of \Cref{thm:PD-TLBT bounds}, we observed that the error bounds for the approximate posterior covariance and posterior mean depend on the prior covariance $\PrCov$ and the observation noise precision $\ObsCov^{-1}$. In our experiment, we study the dependence on $\PrCov$, by providing three experiments with three differently scaled, rank-deficient prior covariances. Let the matrix $\mathbf{P}$ be the unique solution of the reachability Lyapunov equation $\mathbf{A}\mathbf{P} + \mathbf{P}\mathbf{A}^\mathrm{T} = -\mathbf{BB}^\mathrm{T}$ with $\mathbf{B}$ given by the ISS1R benchmark. The solution $\mathbf{P}$ exists because the system matrix $\mbf{A}$ is stable \cite[Proposition 4.27]{Antoulas2005Book}. For a detailed explanation of why this matrix can be chosen as the prior covariance, we refer to \cite{Qian2021Balancing}. We now proceed as in \cite{2025PD-BT} and draw 90 samples from the Gaussian distribution $\mathcal{N}(\mathbf{0},\mathbf{P})$. The empirical covariance of these samples is an approximation of $\mathbf{P}$ of lower rank $s = 89$ and denoted by $\mbf{\Gamma}$. The empirical covariance matrix $\mbf{\Gamma}$ is only positive \textit{semi}definite and thus singular. The bounds we provide in \cref{sec:Lipschitz bounds,sec:PD_BT_bound} allow for singular prior covariances. For our experiments, we consider the prior covariance to be a multiple $\lambda\cdot\mbf{\Gamma}$ of $\mbf{\Gamma}$ with $\lambda \in \{0.01,1,100\}$.

We compute the posterior mean and covariance approximations $\hatPosMn(\data)$ from \eqref{eq:pos_mean_approx} and $\hatPosCov$ from \eqref{eq:pos_covariance_approx} for the initial condition $\mathbf{p}$ using reduced forward models obtained via PD-BT and PD-TLBT as introduced in \cref{subsec:PD-BT}. We compare the error in the approximations to the respective error bounds in \Cref{fig:ISS PDBT}.
For both PD-BT and PD-TLBT, we plot the absolute error $\norm{\PosMn(\data)- \widehat{\pmb{\mu}}_\pos(\data)}_2$ between the full posterior mean and its approximation in the subfigures in the left column, and we plot the absolute error $\norm{\PosCov - \widehat{\mbf{\Gamma}}_\pos}_F$ between the full posterior covariance and its approximation in the subfigures in the right column. We additionally plot the respective error bounds: For PD-BT, this includes \eqref{eq:Frobenius PosMean bound} for the posterior mean approximation and \eqref{eq:Frobenius PosCov bound} for the posterior covariance approximation. Similarly, for PD-TLBT, we plot the bound \eqref{eq:TL Frobenius PosMean bound} for the posterior mean approximation and \eqref{eq:TL Frobenius PosCov bound} for the posterior covariance approximation. 
\begin{figure}[htb]
    \centering
%
%
\definecolor{mycolor1}{rgb}{0.86275,0.14902,0.49804}%
\definecolor{mycolor2}{rgb}{0.39216,0.56078,1.00000}%
\begin{tikzpicture}
\begin{axis}[%
width=0.38\textwidth,
height=1.8in,
at={(0.42\textwidth,4in)},
scale only axis,
xmin=0,
xmax=43,
ymode=log,
ymin=1e-06,
ymax=100000,
yminorticks=true,
xmajorgrids,
ymajorgrids,
xticklabels=\empty,
yticklabels=\empty,
ylabel style={font=\color{white!15!black}},
ylabel={$\mathbf{\Gamma}_\mathrm{pr}=0.01\cdot\mathbf{\Gamma}$},
axis y line*=right,
axis background/.style={fill=white},
title style={font=\bfseries},
title={$\PosCov$ error and bound}
]
\addplot [color=mycolor1, dashdotted, line width=4.0pt, forget plot]
  table[row sep=crcr]{%
1	72.8462404828158\\
3	56.0343348254523\\
5	50.2877325926677\\
7	37.1114286228259\\
9	29.7715850739147\\
11	33.0834487410604\\
13	21.0586320573184\\
15	15.3759895604022\\
17	9.14824896082178\\
19	10.4284375753214\\
21	6.20225069587739\\
23	2.40801831209308\\
25	1.63439968852366\\
27	0.905149455356149\\
31	0.7785615263142\\
33	0.64307343923071\\
35	0.422129880482885\\
37	0.503891018396241\\
39	0.458039724149575\\
41	0.241762302896125\\
43	0.212790352381677\\
};
\addplot [color=mycolor2, dashed, line width=3.0pt, mark options={solid, mycolor2}, forget plot]
  table[row sep=crcr]{%
1	122.348560942444\\
3	85.1152923457292\\
5	69.8758066983765\\
7	47.4656897367163\\
9	35.115387047863\\
11	31.634699305764\\
13	27.6628467124839\\
15	17.8254243823358\\
17	11.0571222988493\\
19	8.87656490747805\\
21	8.05107069863808\\
23	6.20820794138986\\
25	3.77490608587456\\
27	1.40831229110663\\
31	1.12672728208617\\
33	1.00097731363545\\
35	0.950620226013463\\
37	0.875062536649037\\
39	0.752584041862954\\
41	0.69228277537419\\
43	0.58686637159798\\
};
\addplot [color=mycolor1, line width=2.5pt, mark options={solid, rotate=180, mycolor1}, forget plot]
  table[row sep=crcr]{%
1	0.308041801061108\\
3	0.182222148723658\\
5	0.0280809127754185\\
7	0.0230038558966546\\
9	0.0215746677263796\\
11	0.0205413192093359\\
13	0.00685660249684902\\
15	0.00260693800369562\\
17	0.000787764165187501\\
19	0.000518490230751176\\
21	0.00029043767388699\\
23	0.000199450294313074\\
25	0.00012254101495427\\
27	2.23248797247782e-05\\
29	3.64316951177582e-05\\
31	6.57277445581054e-05\\
33	2.51105647161978e-05\\
35	2.95808583582624e-05\\
37	3.64815712139876e-05\\
39	9.45084491238217e-06\\
41	8.8870121458869e-06\\
43	2.00268592018448e-05\\
};
\addplot [color=mycolor2, line width=2.0pt, mark options={solid, mycolor2}, forget plot]
  table[row sep=crcr]{%
1	0.207378648130206\\
3	0.0568168083406118\\
5	0.0362725724186142\\
7	0.0202702617192778\\
9	0.0149755707760837\\
11	0.014078883026598\\
13	0.0112730911421974\\
15	0.00690712941261287\\
17	0.00807194612190271\\
19	0.000812223268828863\\
21	0.000752411938031444\\
23	0.000518412589920835\\
25	0.000232591194709966\\
27	9.90387203501415e-05\\
29	5.52363211046173e-05\\
31	6.15969111880525e-05\\
33	5.13912781096338e-05\\
35	5.18041880240627e-05\\
37	5.39135779404174e-05\\
39	3.10217058297596e-05\\
41	2.17173497511082e-05\\
43	9.66022100189372e-05\\
};
\end{axis}

\begin{axis}[%
width=0.38\textwidth,
height=1.8in,
at={(0in,4in)},
scale only axis,
xmin=0,
xmax=43,
ymode=log,
ymin=1e-06,
ymax=100000,
yminorticks=true,
xmajorgrids,
ymajorgrids,
xticklabels=\empty,
axis y line*=left,
ylabel style={font=\color{white!15!black}},
ylabel={error magnitude},
axis background/.style={fill=white},
title style={font=\bfseries},
title={$\PosMn(\mathbf{m})$ error and bound}
]
\addplot [color=mycolor1, dashdotted, line width=4.0pt, forget plot]
  table[row sep=crcr]{%
1	72.8462404828158\\
3	56.0343348254523\\
5	50.2877325926677\\
7	37.1114286228259\\
9	29.7715850739147\\
11	33.0834487410604\\
13	21.0586320573184\\
15	15.3759895604022\\
17	9.14824896082178\\
19	10.4284375753214\\
21	6.20225069587739\\
23	2.40801831209308\\
25	1.63439968852366\\
27	0.905149455356149\\
31	0.7785615263142\\
33	0.64307343923071\\
35	0.422129880482885\\
37	0.503891018396241\\
39	0.458039724149575\\
41	0.241762302896125\\
43	0.212790352381677\\
};
\addplot [color=mycolor2, dashed, line width=3.0pt, mark options={solid, mycolor2}, forget plot]
  table[row sep=crcr]{%
1	122.348560942444\\
3	85.1152923457292\\
5	69.8758066983765\\
7	47.4656897367163\\
9	35.115387047863\\
11	31.634699305764\\
13	27.6628467124839\\
15	17.8254243823358\\
17	11.0571222988493\\
19	8.87656490747805\\
21	8.05107069863808\\
23	6.20820794138986\\
25	3.77490608587456\\
27	1.40831229110663\\
31	1.12672728208617\\
33	1.00097731363545\\
35	0.950620226013463\\
37	0.875062536649037\\
39	0.752584041862954\\
41	0.69228277537419\\
43	0.58686637159798\\
};
\addplot [color=mycolor1, line width=2.5pt, mark options={solid, rotate=180, mycolor1}, forget plot]
  table[row sep=crcr]{%
1	0.56625308914471\\
3	0.534549678747688\\
5	0.324009135173285\\
7	0.285162768966369\\
9	0.173897926224096\\
11	0.153078884117897\\
13	0.0529962676635618\\
15	0.0539530256068029\\
17	0.045965848204177\\
19	0.0469976742046326\\
21	0.0175650125859476\\
23	0.0110439912188836\\
25	0.0104126563745815\\
27	0.00185577641088562\\
29	0.00139015654874399\\
31	0.00151274230375448\\
33	0.00200095271762843\\
35	0.00150656205554204\\
37	0.00177968006930975\\
39	0.00221541821795821\\
41	0.000907554948713717\\
43	0.000724848237493434\\
};
\addplot [color=mycolor2, line width=2.0pt, mark options={solid, mycolor2}, forget plot]
  table[row sep=crcr]{%
1	0.578738292618722\\
3	0.347357443548839\\
5	0.349801236177076\\
7	0.311200648308104\\
9	0.307235987243932\\
11	0.296813745145706\\
13	0.267221682606829\\
15	0.0536526864445777\\
17	0.0465861906212517\\
19	0.040211727284874\\
21	0.0418625991087682\\
23	0.0314787608631533\\
25	0.0111083827725076\\
27	0.00220427971728969\\
29	0.00190405873861018\\
31	0.00205376169095451\\
33	0.00135186961285087\\
35	0.00114916924632721\\
37	0.00123677360400532\\
39	0.00111107670450636\\
41	0.0010392002029754\\
43	0.00167235496028363\\
};
\end{axis}

\begin{axis}[%
width=0.38\textwidth,
height=1.8in,
at={(0.42\textwidth,2in)},
scale only axis,
xmin=0,
xmax=43,
ymode=log,
ymin=1e-06,
ymax=100000,
yminorticks=true,
xmajorgrids,
ymajorgrids,
xticklabels=\empty,
yticklabels=\empty,
axis y line*=right,
ylabel style={font=\color{white!15!black}},
ylabel={$\mathbf{\Gamma}_\mathrm{pr}=1\cdot\mathbf{\Gamma}$},
axis background/.style={fill=white}
]
\addplot [color=mycolor1, dashdotted, line width=4.0pt, forget plot]
  table[row sep=crcr]{%
1	728.462404828158\\
3	560.343348254523\\
5	502.877325926686\\
7	371.114286228253\\
9	297.715850739137\\
11	330.834487410593\\
13	210.586320573163\\
15	153.759895604037\\
17	91.4824896081714\\
19	104.284375753332\\
21	62.0225069588084\\
23	24.0801831209238\\
25	16.3439968858265\\
27	9.05149455399539\\
29	8.4130114598096\\
31	7.78561526332415\\
33	6.43073439329557\\
35	4.22129879826535\\
37	5.03891018366925\\
39	4.58039724310761\\
41	2.4176230379967\\
43	2.1279035285158\\
};
\addplot [color=mycolor2, dashed, line width=3.0pt, mark options={solid, mycolor2}, forget plot]
  table[row sep=crcr]{%
1	1223.48560942444\\
3	851.152923457292\\
5	698.758066983753\\
7	474.656897367163\\
9	351.15387047863\\
11	316.346993057645\\
13	276.628467124843\\
15	178.254243823355\\
17	110.571222988493\\
19	88.7656490747805\\
21	80.5107069863808\\
23	62.0820794138986\\
25	37.749060858745\\
27	14.0831229110659\\
29	12.6410536210279\\
31	11.2672728208615\\
33	10.0097731363543\\
35	9.50620226013432\\
37	8.75062536649051\\
39	7.52584041862954\\
41	6.9228277537419\\
43	5.8686637159798\\
};
\addplot [color=mycolor1, line width=2.5pt, mark options={solid, rotate=180, mycolor1}, forget plot]
  table[row sep=crcr]{%
1	38.7013681837289\\
3	21.9468536820055\\
5	5.8128807950586\\
7	5.37404357686952\\
9	3.6607114024315\\
11	3.33479971965008\\
13	2.24481786336171\\
15	1.31396541266261\\
17	1.20149421960998\\
19	1.1408733839814\\
21	1.02832769910225\\
23	0.986757445870385\\
25	0.595427665608397\\
27	0.00890042222685355\\
29	0.00717559637242601\\
31	0.0057671367181062\\
33	0.00327304750259325\\
35	0.00317740410802734\\
37	0.00241586446440069\\
39	0.00322831989564457\\
41	0.0023234939135899\\
43	0.00129514529358444\\
};
\addplot [color=mycolor2, line width=2.0pt, mark options={solid, mycolor2}, forget plot]
  table[row sep=crcr]{%
1	26.1036185065742\\
3	7.73182242472281\\
5	5.1006420800371\\
7	3.66156300523649\\
9	3.09757582969444\\
11	2.73850548421672\\
13	2.26904303572294\\
15	1.92490072170926\\
17	0.690664337655903\\
19	0.461118951013079\\
21	0.455732488098935\\
23	0.319138390104768\\
25	0.115330780925029\\
27	0.0160278622706269\\
29	0.00990846779471226\\
31	0.00925298708264972\\
33	0.0100288859846621\\
35	0.0072890057047078\\
37	0.00722417318216092\\
39	0.0134655436317412\\
41	0.00681160115325106\\
43	0.00735215060137307\\
};
\end{axis}

\begin{axis}[%
width=0.38\textwidth,
height=1.8in,
at={(0in,2in)},
scale only axis,
xmin=0,
xmax=43,
ymode=log,
ymin=1e-06,
ymax=100000,
yminorticks=true,
xmajorgrids,
ymajorgrids,
axis y line*=left,
ylabel style={font=\color{white!15!black}},
ylabel={error magnitude},
xticklabels=\empty,
axis background/.style={fill=white}
]
\addplot [color=mycolor1, dashdotted, line width=4.0pt, forget plot]
  table[row sep=crcr]{%
1	728.462404828158\\
3	560.343348254523\\
5	502.877325926686\\
7	371.114286228253\\
9	297.715850739137\\
11	330.834487410593\\
13	210.586320573163\\
15	153.759895604037\\
17	91.4824896081714\\
19	104.284375753332\\
21	62.0225069588084\\
23	24.0801831209238\\
25	16.3439968858265\\
27	9.05149455399539\\
29	8.4130114598096\\
31	7.78561526332415\\
33	6.43073439329557\\
35	4.22129879826535\\
37	5.03891018366925\\
39	4.58039724310761\\
41	2.4176230379967\\
43	2.1279035285158\\
};
\addplot [color=mycolor2, dashed, line width=3.0pt, mark options={solid, mycolor2}, forget plot]
  table[row sep=crcr]{%
1	1223.48560942444\\
3	851.152923457292\\
5	698.758066983753\\
7	474.656897367163\\
9	351.15387047863\\
11	316.346993057645\\
13	276.628467124843\\
15	178.254243823355\\
17	110.571222988493\\
19	88.7656490747805\\
21	80.5107069863808\\
23	62.0820794138986\\
25	37.749060858745\\
27	14.0831229110659\\
29	12.6410536210279\\
31	11.2672728208615\\
33	10.0097731363543\\
35	9.50620226013432\\
37	8.75062536649051\\
39	7.52584041862954\\
41	6.9228277537419\\
43	5.8686637159798\\
};
\addplot [color=mycolor1, line width=2.5pt, mark options={solid, rotate=180, mycolor1}, forget plot]
  table[row sep=crcr]{%
1	4.34307269415546\\
3	3.77560201031165\\
5	2.93350091633671\\
7	2.58166563061304\\
9	1.76773998280689\\
11	1.80031004271934\\
13	1.49227972854382\\
15	1.36099784112444\\
17	1.28479194151258\\
19	1.37216537231034\\
21	1.04107577015117\\
23	1.12437956840519\\
25	1.73853046211641\\
27	0.150883251810569\\
29	0.19856788298962\\
31	0.0925051339923077\\
33	0.0953202690390749\\
35	0.0651709415216885\\
37	0.0579901620450775\\
39	0.0992810929401433\\
41	0.036304583643705\\
43	0.0217291975816447\\
};
\addplot [color=mycolor2, line width=2.0pt, mark options={solid, mycolor2}, forget plot]
  table[row sep=crcr]{%
1	4.414869585251\\
3	3.77660391530314\\
5	3.84177142053376\\
7	2.80006327872224\\
9	2.50835597316991\\
11	3.06512995395833\\
13	1.97199016375573\\
15	1.02065308626817\\
17	1.01065573616377\\
19	1.02073807817868\\
21	0.837051661888445\\
23	0.853555751645615\\
25	0.422185985967014\\
27	0.183770647799172\\
29	0.101217391602939\\
31	0.116943376566102\\
33	0.0909059824517941\\
35	0.0955462614038093\\
37	0.0937565233968122\\
39	0.0962271242266251\\
41	0.0718213758799908\\
43	0.0584687837717382\\
};
\end{axis}

\begin{axis}[%
width=0.38\textwidth,
height=1.8in,
at={(0.42\textwidth,0in)},
scale only axis,
xmin=0,
xmax=43,
xlabel style={font=\color{white!15!black}},
xlabel={$r$},
ymode=log,
ymin=1e-06,
ymax=100000,
yminorticks=true,
xmajorgrids,
ymajorgrids,
yticklabels=\empty,
axis y line*=right,
ylabel style={font=\color{white!15!black}},
ylabel={$\mathbf{\Gamma}_\mathrm{pr}=100\cdot\mathbf{\Gamma}$},
axis background/.style={fill=white},
legend style={font=\footnotesize, at={(0.975in,-0.51in)}, anchor=center, legend cell align=left, align=center, fill=white, draw=black}
]
\addplot [color=mycolor1, dashdotted, line width=4.0pt]
  table[row sep=crcr]{%
1	7284.62404828158\\
3	5603.43348254504\\
5	5028.77325926677\\
7	3711.14286228253\\
9	2977.15850739147\\
11	3308.34487410577\\
13	2105.86320573177\\
15	1537.59895604097\\
17	914.824896082567\\
19	1042.84375753141\\
21	620.225069587485\\
23	240.801831212129\\
25	163.439968861067\\
27	90.5149455556337\\
31	77.8561526360578\\
33	64.3073439262242\\
35	42.2129880514275\\
37	50.389101865044\\
39	45.80397240943\\
41	24.1762304113906\\
43	21.2790352470969\\
};
\addlegendentry{error bound PD-TLBT}

\addplot [color=mycolor2, dashed, line width=3.0pt, mark options={solid, mycolor2}]
  table[row sep=crcr]{%
1	12234.8560942444\\
3	8511.52923457292\\
5	6987.58066983765\\
7	4746.56897367163\\
9	3511.5387047863\\
11	3163.4699305764\\
13	2766.28467124839\\
15	1782.54243823358\\
17	1105.71222988495\\
19	887.656490747805\\
21	805.107069863821\\
23	620.820794138976\\
25	377.49060858745\\
27	140.831229110659\\
31	112.672728208617\\
33	100.097731363543\\
35	95.0620226013432\\
37	87.5062536649037\\
39	75.2584041862954\\
41	69.228277537419\\
43	58.686637159798\\
};
\addlegendentry{error bound PD-BT}

\addplot [color=mycolor1, line width=2.5pt, mark options={solid, rotate=180, mycolor1}]
  table[row sep=crcr]{%
1	3894.15855715617\\
3	2224.23236365179\\
5	663.790450517307\\
7	625.141668787277\\
9	479.209242084975\\
11	447.699168284917\\
13	366.55564269062\\
15	313.203080921195\\
17	307.262763027872\\
19	304.032409532723\\
21	292.640566644919\\
23	285.908581186762\\
25	172.281718341844\\
27	11.8057990903281\\
29	7.59971507266721\\
31	5.89620329547327\\
33	4.29612988046575\\
35	3.0835302607217\\
37	2.04385245379591\\
39	1.6287908849802\\
41	1.35927453297326\\
43	0.909993119097806\\
};
\addlegendentry{$\|\mathbf{\Gamma}_\mathrm{pos}-\hat{\mathbf{\Gamma}}_\mathrm{pos}\|_F$ PD-TLBT}

\addplot [color=mycolor2, line width=2.0pt, mark options={solid, mycolor2}]
  table[row sep=crcr]{%
1	2639.0147112997\\
3	840.882976639554\\
5	596.876441445418\\
7	476.565836447299\\
9	431.333876261585\\
11	339.164208245093\\
13	231.176551751334\\
15	197.947431505832\\
17	79.1931352817822\\
19	53.2015376221748\\
21	47.5805372877771\\
23	27.89119597797\\
25	17.5556446972592\\
27	10.5181421325142\\
29	8.98120314643938\\
31	7.45941358433086\\
33	6.01595560558247\\
35	5.45571149813768\\
37	5.35876080634419\\
39	5.35863358572561\\
41	5.12147678607436\\
43	3.59931110934144\\
};
\addlegendentry{$\|\mathbf{\Gamma}_\mathrm{pos}-\hat{\mathbf{\Gamma}}_\mathrm{pos}\|_F$ PD-BT}

\end{axis}

\begin{axis}[%
width=0.38\textwidth,
height=1.8in,
at={(0in,0in)},
scale only axis,
xmin=0,
xmax=43,
xlabel style={font=\color{white!15!black}},
xlabel={$r$},
ymode=log,
ymin=1e-06,
ymax=100000,
yminorticks=true,
xmajorgrids,
ymajorgrids,
axis y line*=left,
ylabel style={font=\color{white!15!black}},
ylabel={error magnitude},
axis background/.style={fill=white},
legend style={font=\footnotesize, at={(0.95in,-0.525in)}, anchor=center, legend cell align=left, align=center, fill=white, draw=black}
]
\addplot [color=mycolor1, dashdotted, line width=4.0pt,mark options={solid, mycolor1}]
  table[row sep=crcr]{%
1	7284.62404828158\\
3	5603.43348254504\\
5	5028.77325926677\\
7	3711.14286228253\\
9	2977.15850739147\\
11	3308.34487410577\\
13	2105.86320573177\\
15	1537.59895604097\\
17	914.824896082567\\
19	1042.84375753141\\
21	620.225069587485\\
23	240.801831212129\\
25	163.439968861067\\
27	90.5149455556337\\
31	77.8561526360578\\
33	64.3073439262242\\
35	42.2129880514275\\
37	50.389101865044\\
39	45.80397240943\\
41	24.1762304113906\\
43	21.2790352470969\\
};
\addlegendentry{error bound PD-TLBT}

\addplot [color=mycolor2, dashed, line width=3.0pt, mark options={solid, mycolor2}]
  table[row sep=crcr]{%
1	12234.8560942444\\
3	8511.52923457292\\
5	6987.58066983765\\
7	4746.56897367163\\
9	3511.5387047863\\
11	3163.4699305764\\
13	2766.28467124839\\
15	1782.54243823358\\
17	1105.71222988495\\
19	887.656490747805\\
21	805.107069863821\\
23	620.820794138976\\
25	377.49060858745\\
27	140.831229110659\\
31	112.672728208617\\
33	100.097731363543\\
35	95.0620226013432\\
37	87.5062536649037\\
39	75.2584041862954\\
41	69.228277537419\\
43	58.686637159798\\
};
\addlegendentry{error bound PD-BT}

\addplot [color=mycolor1, line width=2.5pt, mark options={solid, rotate=180, mycolor1}]
  table[row sep=crcr]{%
1	104.213322123559\\
3	99.2864592020448\\
5	67.3957056907097\\
7	26.5891400181545\\
9	25.674506048869\\
11	32.8995511888322\\
13	23.6270896791511\\
15	18.7080417868572\\
17	17.6751652197013\\
19	17.7504959479493\\
21	13.5209757722367\\
23	18.1561530543059\\
25	25.5088971744777\\
27	4.20397551505067\\
29	2.98566811489581\\
31	3.77028528365213\\
33	2.26096204851783\\
35	1.61457187406317\\
37	2.99999535627274\\
39	1.71141629009551\\
41	1.52937213277746\\
43	1.50609053774351\\
};
\addlegendentry{$\|\pmb{\mu}_\mathrm{pos}-\hat{\pmb{\mu}}_\mathrm{pos}\|_2$ PD-TLBT}

\addplot [color=mycolor2, line width=2.0pt, mark options={solid, mycolor2}]
  table[row sep=crcr]{%
1	103.725658053289\\
3	77.6948361538392\\
5	76.8826735675827\\
7	65.4716089180722\\
9	27.318981678299\\
11	24.920395089763\\
13	18.4767503634088\\
15	20.4863714061125\\
17	13.8612619119744\\
19	10.3546165391561\\
21	11.5005539004305\\
23	9.17099940789314\\
25	7.0854112827495\\
27	2.76881551373891\\
29	2.66179646483131\\
31	4.11856121653644\\
33	2.69137991825626\\
35	2.78672808234154\\
37	4.39884427563836\\
39	2.97881355979589\\
41	3.0188025632457\\
43	3.79464379417888\\
};
\addlegendentry{$\|\pmb{\mu}_\mathrm{pos}-\hat{\pmb{\mu}}_\mathrm{pos}\|_2$ PD-BT}

\end{axis}
\end{tikzpicture}

    \caption{\small ISS1R benchmark: Comparison of the actual mean and covariance approximation errors by PD-BT and PD-TLBT with the respective error bounds for different prior covariances. Note that the curves of the error bounds have been scaled to facilitate comparison with the curves of the actual errors; see \Cref{remark_figures}.}
    \label{fig:ISS PDBT}
\end{figure}

The plots of the actual errors in the posterior mean and covariance approximations and the error bounds are given for the three different prior covariances described above. The plots in the top row show the results for the prior covariance $\PrCov=0.01\cdot\mathbf{\Gamma}$ with the smallest norm, which represents a very certain prior belief. The plots in the middle row show the results for the prior covariance $\PrCov=\mathbf{\Gamma}$ and the plots in the bottom row show the results for the prior covariance $\PrCov=100\cdot\mathbf{\Gamma}$ with the largest norm, which represents a very uncertain prior belief. The actual errors and error bounds are plotted with respect to the rank $r$ of the reduced forward model obtained by PD-(TL)BT. A smaller rank $r$ corresponds to a greater number of truncated HSVs and lower model accuracy.

\paragraph{Two-dimensional heat equation with Gaussian covariance function}
For a second example, we will use the same LTI system in \cite{2025PD-BT} considering the heat equation on the unit square with homogeneous Dirichlet boundary conditions, over the time interval $[0,5]$. One-dimensional versions of this example were studied in \cite{Koe2022TLBTDA,Qian2021Balancing}. We chose a uniform grid discretization of the unit square using points with coordinates $(i\Delta , j\Delta)$ for $\Delta= 1/51$ and $i,j=1,\ldots,50$, and used the standard centered finite differences stencil to discretize the Laplacian with diffusivity constant $\alpha=0.1$.
These choices result in a system matrix $\mbf{A}\in \R^{2500\times 2500}$ and an output matrix $\mbf{C}\in\R^{1\times 2500}$ in \eqref{eq:BayesianSystem_nonoise}, where $\mbf{C}$ is chosen such that the output $\mbf{y}(t)$ is the average temperature over the spatial domain at time $t$.

We construct three prior covariances $\PrCov$ by evaluating the Gaussian covariance function $f(\mbf{z},\mbf{z}')=\sigma_f^2\exp(-\frac{\norm{\mbf{z}-\mbf{z}'}^2}{2\ell^2})$ with signal variance $\sigma_f=1$ on the uniform grid described above and three different length scales $\ell=0.01$,  $\ell=0.1$, and $\ell=1$. 
Given the choice of the Gaussian covariance function, it holds that as $\ell$ increases, the covariance between two distinct points in the spatial domain increases. 
The numerical rank of the prior covariance matrix $\PrCov$ corresponding to the length scales is $s=2500$, $s=609$ and $s=37$, respectively. 
To generate the vector $\data$ of noisy output measurements, we initialized the system at a reference initial condition $\mbf{p}$ drawn from the prior and then simulated the unforced linear system $\mbf{x}(t)$ according to \eqref{eq:BayesianSystem_nonoise} over the time interval $[0,5]$. We computed the observations by calculating matrix exponentials to directly solve the heat equation at the different observation times. Measurements were computed according to \eqref{eq:BayesianSystem_noise}, using the resulting trajectory, equidistant times $t_k=k\Delta t$ for $\Delta t=0.2$, $k=1,\ldots, n\coloneqq 5/ \Delta t$, and $\mbf{\Gamma}_\epsilon=0.08$.

For the computations of the posterior mean and covariance approximations using PD-BT and PD-TLBT, we proceeded in the analogous way as for the ISS1R benchmark example. We compare the error in the approximations to the respective error bounds in \Cref{fig:heat PDBT}.

\begin{figure}[htb]
    \centering
%
%
\definecolor{mycolor1}{rgb}{0.86275,0.14902,0.49804}%
\definecolor{mycolor2}{rgb}{0.39216,0.56078,1.00000}%
\begin{tikzpicture}
\begin{axis}[%
width=0.38\textwidth,
height=1.8in,
at={(0.42\textwidth,4in)},
scale only axis,
xmin=0,
xmax=31,
ymode=log,
ymin=1e-14,
ymax=10000,
yminorticks=true,
xmajorgrids,
ymajorgrids,
xticklabels=\empty,
yticklabels=\empty,
ylabel style={font=\color{white!15!black}},
ylabel={$\ell=0.01$},
axis y line*=right,
axis background/.style={fill=white},
title style={font=\bfseries},
title={$\PosCov$ error and bound},
legend style={font=\footnotesize, at={(0.585,0.84)}, anchor=center, legend cell align=left, align=center, fill=white, draw=black},
]
\addplot [color=mycolor1, dashed, line width=3.0pt]
  table[row sep=crcr]{%
1	0.0374331641716124\\
3	0.00383640838171319\\
5	0.000606625239892379\\
7	0.000125987023651695\\
9	1.67596559621806e-05\\
11	2.90922394056054e-06\\
13	3.74449148389438e-07\\
15	3.71016148596468e-08\\
17	4.8634399040709e-08\\
19	5.09701647854178e-08\\
21	2.49079868423351e-08\\
23	1.8154659134387e-08\\
25	6.28895840297055e-08\\
27	9.68693352105626e-08\\
29	4.5118877665176e-08\\
31	2.20157580024006e-08\\
};
\addlegendentry{error bound PD-TLBT}

\addplot [color=mycolor2, dashed, line width=2.0pt]
  table[row sep=crcr]{%
1	0.0374749300821857\\
3	0.00393808522634834\\
5	0.000657402007973438\\
7	0.000140698524707866\\
9	1.84414189052637e-05\\
11	3.2292298726055e-06\\
13	4.06549847060098e-07\\
15	7.21930075278261e-08\\
17	7.5887497263248e-09\\
19	1.14334388282482e-09\\
21	1.36791307153911e-10\\
23	1.7065330054757e-11\\
25	2.6488840148962e-12\\
27	3.90183847985137e-13\\
29	5.66860050615322e-14\\
31	3.09686693032246e-14\\
};
\addlegendentry{error bound PD-BT}

\addplot [color=mycolor1, line width=3.0pt]
  table[row sep=crcr]{%
1	0.0145082515230466\\
3	0.000213878405309958\\
5	6.20377074058868e-05\\
7	5.97078149790451e-06\\
9	1.73655806304953e-06\\
11	3.65401977047373e-07\\
13	1.59941218570414e-08\\
15	7.90259875774839e-09\\
17	7.90259843028404e-09\\
19	3.46449734802499e-09\\
21	3.46448121042548e-09\\
23	1.36383705054169e-09\\
31	1.36371974645001e-09\\
};
\addlegendentry{$\|\mathbf{\Gamma}_\mathrm{pos}-\hat{\mathbf{\Gamma}}_\mathrm{pos}\|_F$ PD-TLBT}

\addplot [color=mycolor2, line width=2.0pt]
  table[row sep=crcr]{%
1	0.0145082514665207\\
3	0.00021387840986393\\
5	6.20377080759054e-05\\
7	5.97078151166252e-06\\
9	1.73655805060215e-06\\
11	3.65401913673608e-07\\
13	1.59940207252749e-08\\
15	7.90250482396454e-09\\
17	1.36449410457239e-09\\
19	1.47917719520613e-10\\
21	6.238547470906e-12\\
23	3.09765855893512e-12\\
25	3.75588662031667e-13\\
27	5.42616711295398e-14\\
29	4.27379313151743e-14\\
31	4.18975486809579e-14\\
};
\addlegendentry{$\|\mathbf{\Gamma}_\mathrm{pos}-\hat{\mathbf{\Gamma}}_\mathrm{pos}\|_F$ PD-BT}

\end{axis}

\begin{axis}[%
width=0.38\textwidth,
height=1.8in,
at={(0in,4in)},
scale only axis,
xmin=0,
xmax=31,
ymode=log,
ymin=1e-14,
ymax=10000,
yminorticks=true,
xmajorgrids,
ymajorgrids,
xticklabels=\empty,
axis y line*=left,
ylabel style={font=\color{white!15!black}},
ylabel={error magnitude},
axis background/.style={fill=white},
title style={font=\bfseries},
title={$\PosMn(\mathbf{m})$ error and bound},
legend style={font=\footnotesize, at={(0.585,0.85)}, anchor=center, legend cell align=left, align=center, fill=white, draw=black}
]
\addplot [color=mycolor1, dashed, line width=3.0pt]
  table[row sep=crcr]{%
1	0.0374331641716124\\
3	0.00383640838171319\\
5	0.000606625239892379\\
7	0.000125987023651695\\
9	1.67596559621806e-05\\
11	2.90922394056054e-06\\
13	3.74449148389438e-07\\
15	3.71016148596468e-08\\
17	4.8634399040709e-08\\
19	5.09701647854178e-08\\
21	2.49079868423351e-08\\
23	1.8154659134387e-08\\
25	6.28895840297055e-08\\
27	9.68693352105626e-08\\
29	4.5118877665176e-08\\
31	2.20157580024006e-08\\
};
\addlegendentry{error bound PD-TLBT}

\addplot [color=mycolor2, dashed, line width=2.0pt]
  table[row sep=crcr]{%
1	0.0374749300821857\\
3	0.00393808522634834\\
5	0.000657402007973438\\
7	0.000140698524707866\\
9	1.84414189052637e-05\\
11	3.2292298726055e-06\\
13	4.06549847060098e-07\\
15	7.21930075278261e-08\\
17	7.5887497263248e-09\\
19	1.14334388282482e-09\\
21	1.36791307153911e-10\\
23	1.7065330054757e-11\\
25	2.6488840148962e-12\\
27	3.90183847985137e-13\\
29	5.66860050615322e-14\\
31	3.09686693032246e-14\\
};
\addlegendentry{error bound PD-BT}

\addplot [color=mycolor1, line width=3.0pt]
  table[row sep=crcr]{%
1	0.0343008419059636\\
3	0.00260513910488139\\
5	0.000503021426448895\\
7	5.41607123864395e-05\\
9	2.65189343713842e-06\\
11	8.75027698539655e-07\\
13	1.0616272747482e-07\\
15	5.55970125639415e-08\\
17	5.55970099220346e-08\\
19	2.23936902831793e-08\\
21	2.23935817519642e-08\\
23	8.62706987238125e-09\\
31	8.62627164558365e-09\\
};
\addlegendentry{$\|\pmb{\mu}_\mathrm{pos}-\hat{\pmb{\mu}}_\mathrm{pos}\|_2$ PD-TLBT}

\addplot [color=mycolor2, line width=2.0pt]
  table[row sep=crcr]{%
1	0.0343008417299671\\
3	0.00260513911775241\\
5	0.000503021427096988\\
7	5.41607121133964e-05\\
9	2.65189348801498e-06\\
11	8.75027598307593e-07\\
13	1.0616271439343e-07\\
15	5.5595722410025e-08\\
17	8.63113947516503e-09\\
19	1.00913585437134e-09\\
21	5.33519140932417e-11\\
23	6.42100267293979e-12\\
25	8.36014207737623e-13\\
27	3.50477673652056e-13\\
29	1.10532430450812e-13\\
31	1.0033727770259e-13\\
};
\addlegendentry{$\|\pmb{\mu}_\mathrm{pos}-\hat{\pmb{\mu}}_\mathrm{pos}\|_2$ PD-BT}

\end{axis}

\begin{axis}[%
width=0.38\textwidth,
height=1.8in,
at={(0.42\textwidth,2in)},
scale only axis,
xmin=0,
xmax=31,
ymode=log,
ymin=1e-14,
ymax=10000,
yminorticks=true,
xmajorgrids,
ymajorgrids,
xticklabels=\empty,
yticklabels=\empty,
axis y line*=right,
ylabel style={font=\color{white!15!black}},
ylabel={$\ell=0.1$},
axis background/.style={fill=white}
]
\addplot [color=mycolor1, dashed, line width=3.0pt, forget plot]
  table[row sep=crcr]{%
1	26.0440090483331\\
3	1.12705556324678\\
5	0.108147874357833\\
7	0.00823660873707546\\
9	0.000782144027792437\\
11	5.36625849884456e-05\\
13	3.98856906166877e-05\\
15	2.63580403271126e-05\\
17	2.0807638809517e-05\\
19	4.14941266265397e-05\\
21	3.02697803048571e-05\\
23	1.82282737725748e-05\\
25	3.17313517044856e-06\\
27	2.35326002504305e-05\\
29	2.22119461931399e-05\\
31	8.39532653739943e-06\\
};
\addplot [color=mycolor2, dashed, line width=2.0pt, forget plot]
  table[row sep=crcr]{%
1	26.0799944531135\\
3	1.1724889129518\\
5	0.125679618593549\\
7	0.00948070881280278\\
9	0.000879648463705659\\
11	6.48879994057417e-05\\
13	5.53124432817388e-06\\
15	3.18074985652974e-07\\
17	3.26857900944849e-08\\
19	3.80259077924966e-09\\
21	4.51400678420974e-10\\
23	3.97181058176804e-11\\
25	2.05768365511011e-11\\
27	2.52965545920979e-11\\
29	2.32600155140991e-11\\
31	6.70416376249006e-11\\
};
\addplot [color=mycolor1, line width=3.0pt, forget plot]
  table[row sep=crcr]{%
1	13.7338358493693\\
3	0.125190708077601\\
5	0.00959022970369609\\
7	0.000693359553585816\\
9	4.63857530241805e-05\\
11	1.23835235102388e-06\\
13	8.54547743313206e-07\\
15	3.52906701837321e-07\\
17	3.34597283091326e-07\\
19	3.39279835728181e-07\\
21	3.40010658315284e-07\\
23	1.07916388429362e-07\\
29	1.07481045346579e-07\\
31	1.07189718251911e-07\\
};
\addplot [color=mycolor2, line width=2.0pt, forget plot]
  table[row sep=crcr]{%
1	13.7338357941125\\
3	0.125190709346021\\
5	0.009590229795767\\
7	0.000693359547505332\\
9	4.63857632550417e-05\\
11	1.23700126866081e-06\\
13	3.35317861673538e-07\\
15	1.98098450517376e-08\\
17	1.10256234994963e-09\\
19	3.7025079868858e-10\\
21	4.88288202697927e-11\\
23	4.55419634288214e-11\\
25	4.53495955386323e-11\\
27	3.16549684228388e-11\\
29	2.28672491753071e-11\\
31	2.42307385270627e-11\\
};
\end{axis}

\begin{axis}[%
width=0.38\textwidth,
height=1.8in,
at={(0in,2in)},
scale only axis,
xmin=0,
xmax=31,
ymode=log,
ymin=1e-14,
ymax=10000,
yminorticks=true,
xmajorgrids,
ymajorgrids,
axis y line*=left,
ylabel style={font=\color{white!15!black}},
ylabel={error magnitude},
xticklabels=\empty,
axis background/.style={fill=white}
]
\addplot [color=mycolor1, dashed, line width=3.0pt, forget plot]
  table[row sep=crcr]{%
1	26.0440090483331\\
3	1.12705556324678\\
5	0.108147874357833\\
7	0.00823660873707546\\
9	0.000782144027792437\\
11	5.36625849884456e-05\\
13	3.98856906166877e-05\\
15	2.63580403271126e-05\\
17	2.0807638809517e-05\\
19	4.14941266265397e-05\\
21	3.02697803048571e-05\\
23	1.82282737725748e-05\\
25	3.17313517044856e-06\\
27	2.35326002504305e-05\\
29	2.22119461931399e-05\\
31	8.39532653739943e-06\\
};
\addplot [color=mycolor2, dashed, line width=2.0pt, forget plot]
  table[row sep=crcr]{%
1	26.0799944531135\\
3	1.1724889129518\\
5	0.125679618593549\\
7	0.00948070881280278\\
9	0.000879648463705659\\
11	6.48879994057417e-05\\
13	5.53124432817388e-06\\
15	3.18074985652974e-07\\
17	3.26857900944849e-08\\
19	3.80259077924966e-09\\
21	4.51400678420974e-10\\
23	3.97181058176804e-11\\
25	2.05768365511011e-11\\
27	2.52965545920979e-11\\
29	2.32600155140991e-11\\
31	6.70416376249006e-11\\
};
\addplot [color=mycolor1, line width=3.0pt, forget plot]
  table[row sep=crcr]{%
1	0.726629897136858\\
3	0.020881575026148\\
5	0.000934586512617077\\
7	0.000240345393259064\\
9	1.65198739354395e-05\\
11	9.32535086395869e-07\\
13	2.54610455387159e-07\\
15	3.13173299016707e-08\\
17	2.79972643393238e-08\\
21	2.82294455173293e-08\\
23	5.37053563377628e-09\\
25	5.1983809910359e-09\\
29	5.10863244436117e-09\\
31	5.07968351590653e-09\\
};
\addplot [color=mycolor2, line width=2.0pt, forget plot]
  table[row sep=crcr]{%
1	0.726629896079157\\
3	0.0208815749185139\\
5	0.000934586525642424\\
7	0.000240345393742877\\
9	1.65198731843735e-05\\
11	9.32563216067913e-07\\
13	2.69952190659782e-08\\
15	3.97293815265182e-10\\
17	5.71204652050001e-11\\
19	1.39299586538611e-11\\
21	4.30870047070895e-12\\
23	1.38488331982158e-12\\
25	1.37364985627454e-12\\
27	1.02178066004127e-12\\
29	6.0267900165765e-13\\
31	1.58201166230786e-12\\
};
\end{axis}

\begin{axis}[%
width=0.38\textwidth,
height=1.8in,
at={(0.42\textwidth,0in)},
scale only axis,
xmin=0,
xmax=31,
xlabel style={font=\color{white!15!black}},
xlabel={$r$},
ymode=log,
ymin=1e-14,
ymax=10000,
yminorticks=true,
xmajorgrids,
ymajorgrids,
yticklabels=\empty,
axis y line*=right,
ylabel style={font=\color{white!15!black}},
ylabel={$\ell=1$},
axis background/.style={fill=white}
]
\addplot [color=mycolor1, dashed, line width=3.0pt, forget plot]
  table[row sep=crcr]{%
1	595.976487488812\\
3	15.0467859964282\\
5	0.425408135531475\\
7	0.0314011980312776\\
9	0.00048015764248206\\
11	0.0011333580384808\\
13	0.000702878817487556\\
15	0.000725930121208339\\
17	0.00048015764248206\\
19	0.000654344568605116\\
21	0.00104253776457257\\
23	0.00117614122017876\\
25	0.000851228520261534\\
27	0.00160281030904412\\
29	0.000362965060604169\\
31	0.000925382963386527\\
};
\addplot [color=mycolor2, dashed, line width=2.0pt, forget plot]
  table[row sep=crcr]{%
1	599.734661182398\\
3	16.2499305879049\\
5	0.48994437462502\\
7	0.0381686975544854\\
9	0.00092634910203109\\
11	5.25819922693556e-05\\
13	3.44521985188038e-06\\
15	1.0618628130594e-07\\
17	7.26506834531048e-09\\
19	4.14289209349607e-10\\
21	1.5885153529767e-10\\
23	1.21818230974305e-10\\
25	7.88233041188475e-11\\
27	2.48839498204795e-11\\
29	6.35413240500584e-11\\
31	3.2148887459358e-10\\
};
\addplot [color=mycolor1, line width=3.0pt, forget plot]
  table[row sep=crcr]{%
1	2.82305652877869\\
3	0.0227891282030315\\
5	0.00135490523509789\\
7	7.6421568150468e-05\\
9	4.50334372395982e-06\\
11	2.5697002849755e-07\\
13	2.53390666926685e-07\\
15	8.6419284127184e-08\\
17	8.65409453787879e-08\\
19	9.96365204400632e-08\\
21	1.82866134426189e-08\\
23	1.84376268116695e-08\\
25	1.80438490033423e-08\\
27	3.93970159446672e-09\\
29	3.92183011921686e-09\\
31	5.12537018720632e-09\\
};
\addplot [color=mycolor2, line width=2.0pt, forget plot]
  table[row sep=crcr]{%
1	2.82305653532384\\
3	0.0227891281444658\\
5	0.00135490523836262\\
7	7.64215696476966e-05\\
9	4.50333726585331e-06\\
11	2.53052157478075e-07\\
13	1.79041177571703e-08\\
15	5.61521461144101e-10\\
17	6.65695747816781e-11\\
19	4.98446093835021e-11\\
21	4.93668446313101e-11\\
23	6.66373411349337e-12\\
25	6.56483223492918e-12\\
29	6.15598207114926e-12\\
31	6.14871605776118e-12\\
};
\end{axis}

\begin{axis}[%
width=0.38\textwidth,
height=1.8in,
at={(0in,0in)},
scale only axis,
xmin=0,
xmax=31,
xlabel style={font=\color{white!15!black}},
xlabel={$r$},
ymode=log,
ymin=1e-14,
ymax=10000,
yminorticks=true,
xmajorgrids,
ymajorgrids,
axis y line*=left,
ylabel style={font=\color{white!15!black}},
ylabel={error magnitude},
axis background/.style={fill=white}
]
\addplot [color=mycolor1, dashed, line width=3.0pt, forget plot]
  table[row sep=crcr]{%
1	595.976487488812\\
3	15.0467859964282\\
5	0.425408135531475\\
7	0.0314011980312776\\
9	0.00048015764248206\\
11	0.0011333580384808\\
13	0.000702878817487556\\
15	0.000725930121208339\\
17	0.00048015764248206\\
19	0.000654344568605116\\
21	0.00104253776457257\\
23	0.00117614122017876\\
25	0.000851228520261534\\
27	0.00160281030904412\\
29	0.000362965060604169\\
31	0.000925382963386527\\
};
\addplot [color=mycolor2, dashed, line width=2.0pt, forget plot]
  table[row sep=crcr]{%
1	599.734661182398\\
3	16.2499305879049\\
5	0.48994437462502\\
7	0.0381686975544854\\
9	0.00092634910203109\\
11	5.25819922693556e-05\\
13	3.44521985188038e-06\\
15	1.0618628130594e-07\\
17	7.26506834531048e-09\\
19	4.14289209349607e-10\\
21	1.5885153529767e-10\\
23	1.21818230974305e-10\\
25	7.88233041188475e-11\\
27	2.48839498204795e-11\\
29	6.35413240500584e-11\\
31	3.2148887459358e-10\\
};
\addplot [color=mycolor1, line width=3.0pt, forget plot]
  table[row sep=crcr]{%
1	1.96455053127078\\
3	0.0194252330286203\\
5	0.000755165100463845\\
7	6.96817980909181e-05\\
9	2.20033511450863e-06\\
11	4.0139984562925e-08\\
13	4.02300234647599e-08\\
15	2.55712530751161e-08\\
17	2.55687910463276e-08\\
19	2.51130531302904e-08\\
21	4.51411561014435e-09\\
23	4.53503368594126e-09\\
25	4.60397726903076e-09\\
27	4.31315489204748e-10\\
29	4.46268389904074e-10\\
31	6.42579442115763e-10\\
};
\addplot [color=mycolor2, line width=2.0pt, forget plot]
  table[row sep=crcr]{%
1	1.96455054339491\\
3	0.0194252328895781\\
5	0.000755165101008608\\
7	6.96817979722254e-05\\
9	2.20033777117826e-06\\
11	4.02194607138116e-08\\
13	4.55108879708748e-09\\
15	3.30874829464766e-11\\
17	2.34514403826911e-11\\
19	8.00040149352123e-12\\
21	6.67127519602744e-12\\
23	7.0246962383083e-12\\
25	8.20223219133085e-12\\
27	7.58170458452453e-12\\
29	7.08115455589838e-12\\
31	7.15420034711143e-12\\
};
\end{axis}
\end{tikzpicture}

    \caption{\small Two-dimensional heat equation: Comparison of the actual mean and covariance approximation errors by PD-BT and PD-TLBT with the respective error bounds for different prior covariances. Note that the curves of the error bounds have been scaled to facilitate comparison of the curves of the actual errors; see \Cref{remark_figures}.}
    \label{fig:heat PDBT}
\end{figure}

Plots of the actual errors in the posterior mean and covariance approximations, as well as the error bounds, are provided for the three different prior covariances corresponding to the three different choices of length scale $\ell$ described above. The plots in the top, middle and bottom rows of \Cref{fig:heat PDBT} show the results for the prior covariance matrix $\PrCov$ arising from the Gaussian covariance function with the length scale $\ell=0.01$, $\ell=0.1$, and $\ell=1$ respectively. The actual errors and error bounds are plotted with respect to the rank $r$ of the reduced forward model obtained by PD-(TL)BT. Note that a smaller rank $r$ corresponds to a greater number of truncated HSVs and lower model accuracy.

\begin{remark}
 \label{remark_figures}
We make two remarks about the results shown in \Cref{fig:ISS PDBT,fig:heat PDBT}. First, while we plot the actual errors exactly without any scaling, we replace the constants $\kappa$, $C$ and $C'$ in the error bounds with the squared norm term $\| \mbf{L}_\pr\|_\infty^2$ from the definition \eqref{eq:error_covariance_for_perturbed_forward_model_by_ppH_constant} of $C$. This scales down the magnitude of the error bounds in order to bring the curves for the error bounds and the curves for the actual errors closer together, thus enhancing readability of the plots.
Second, the results shown are for one realization of the reference initial condition and the observation noise used to generate the data $\data$, as well as the empirical prior covariance for the ISS1R benchmark example. However, the results of other numerical experiments that we do not show here demonstrate similar behavior for different realizations of these quantities. This is consistent with \Cref{theorem03}, which treats $\PrCov$ and $\data$ as given inputs to the inverse problem, and thus also for \Cref{thm:PD-BT bounds} and \Cref{thm:PD-TLBT bounds}. In particular, varying $\PrCov$ and $\data$ does not change the key implication of the error bounds, namely that the errors in the posterior mean and covariance are of first order in the error of the root prior-preconditioned Hessian.
\end{remark}

The complete code for our experiments is available under \url{https://github.com/joskoUP/PD-BT-bounds}.

\subsection{Results and discussion}\label{subsec:Numerics discussion}
\paragraph{ISS1R benchmark} We observe that for all choices of prior covariance, as well as for both mean and covariance approximations, the decreasing trend in the actual errors matches the decreasing trend in the error bounds from \cref{thm:PD-BT bounds,thm:PD-TLBT bounds}. As the magnitude of the prior covariance increases, i.e., as the factor $\lambda$ increases, so do the actual errors and the error bounds. The actual errors and error bounds for PD-BT and PD-TLBT do not differ significantly in either plot. This is due to the relatively high end time, $t_n = 8$, of the chosen benchmark experiments from \cite{2025PD-BT}. Previous work \cite{Koe2022TLBTDA} has shown that, for smaller end times, the relative errors of PD-TLBT are smaller than the relative errors of PD-BT, and approach the relative errors of PD-BT as the end time $t_n$ increases to infinity.

These findings highlight how the bounds \eqref{eq:Frobenius PosCov bound}--\eqref{eq:TL Frobenius PosMean bound} capture the behavior of the posterior approximation. A priori bounds tend to be conservative, resulting in them being much larger than the actual approximation error. This can also be observed in \cref{fig:ISS PDBT}, where the error bounds are larger than the actual errors, even when the Lipschitz constants $C$ and $C'$ are ignored. The constant $C$ is on the order of $10^{3}$, $10^{6}$, and $10^{9}$ for the three experiments and increases as the magnitude of the prior covariance increases.  Similarly, the constant $C'$ is on the order of $10^{5}$, $10^{9}$, and $10^{14}$, respectively for this example. Furthermore, the plots highlight the importance of prior knowledge in choosing the reduction rank. For the `small' prior covariance $\PrCov=0.01\cdot\mathbf{\Gamma}$ in the top row of \cref{fig:ISS PDBT}, a modest rank of about $r=20$ suffices to achieve an absolute error of around $10^{-2}$ in the posterior mean approximation. In the middle row with $\PrCov=\mathbf{\Gamma}$, a rank of about $r=40$ is required to achieve a similar level of approximation quality. For the `largest' prior covariance $\PrCov=100\cdot\mathbf{\Gamma}$ in the bottom row, the overall approximation quality is much worse and does not even attain an absolute error of $10^{-2}$ for the ranks $r$ that are plotted.

\paragraph{Two-dimensional heat equation with Gaussian covariance function}
The second experiment shows that for all choices of prior covariance and for both mean and covariance approximations, the trend of decreasing actual errors aligns with the trend of decreasing error bounds from \cref{thm:PD-BT bounds,thm:PD-TLBT bounds}. As \cref{fig:heat PDBT} shows, as the length scale $\ell$ increases, so do the actual errors and the error bounds. For this example, the constant $C$ is on the order of $10^{0}$, $10^{3}$, and $10^{6}$ for the three experiments and increases as the length scale increases.  Similarly, the constant $C'$ is on the order of $10^{1}$, $10^{2}$, and $10^{6}$, respectively. After a certain value of the rank $r$, the PD-TLBT errors and bounds appear to approximately plateau, whereas the PD-BT errors and bounds continue to decrease. This is due to the less stable numerical implementation of PD-TLBT, especially when computing the error bound values. This computation is more demanding than that of PD-BT because it requires a matrix exponential rather than just solutions to Sylvester equations. Since the HSVs decay quickly for the heat equation, the posterior errors and bounds decay quickly as well, resulting in numerical inaccuracies for small error magnitudes.

In the preceding two examples, the dependence of the errors and the error bounds on the prior highlights the impact of prior selection on posterior computation and approximation quality. From a system-theoretic point of view, it is clear that the decay of the HSVs determines the quality of the posterior approximation. However, the experiments indicate that choosing a good prior---and possibly prior hyperparameters, as in the case of a Gaussian covariance function---affects the quality of the posterior approximation as well.
This is consistent with the appearance of prior-dependent quantities in the error bounds of \Cref{theorem03}, \Cref{thm:PD-BT bounds}, and \Cref{thm:PD-TLBT bounds}.

The experiments reported above represent the first published implementation of PD-TLBT. While the concept of PD-TLBT was mentioned in \cite{2025PD-BT}, the analysis of Gramians for PD-TLBT has not been studied in the literature, to the best of our knowledge. Time-limited balanced truncation for LTI systems \cite{Wodek1990TLGramians,Kuerschner2018TLBT,Redmann2020TLBT,Redmann2018H2TLBT} was applied to data assimilation \cite{Koe2022TLBTDA,koenig2022TLBT4DVAR} but never to the prior-driven system \eqref{eq:PriorDrivenSys}. We address this gap by introducing PD-TLBT in \cref{subsec:PD-BT} and by providing the corresponding code and experiments in \cref{sec:Numerics}. In the experiments above, we nonetheless use a stable system to compare the two proposed bounds for PD-BT and PD-TLBT, because PD-TLBT can be applied to unstable systems, but PD-BT cannot. However, note that the expressions under the square roots in the bounds \eqref{eq:TL Frobenius PosCov bound} and \eqref{eq:TL Frobenius PosMean bound} for PD-TLBT may become too large to be useful for unstable systems at higher end times due to the matrix exponential in the reachability Gramian $\mathbf{P}_{\mathcal{T}}$ \eqref{eq:TL_Gramians} and its truncated and mixed versions \cite[Section 3.1]{Kuerschner2018TLBT}.

\section{Discussion of further system-theoretic MOR methods in Bayesian inference}\label{sec:MOR Alternatives}
In this section, we discuss whether posterior error bounds analogous to \Cref{thm:PD-BT bounds} and \Cref{thm:PD-TLBT bounds} can be proven for MOR methods other than PD-(TL)BT. 

\paragraph{Balanced truncation based on optimal low-rank approximation} In \cite{Qian2021Balancing}, a balanced truncation-based MOR approach is applied to solve the same Bayesian inverse problem consisting of inferring the initial condition of the system \eqref{eq:BayesianSystem}, and thus shares a common goal with PD-BT.
This approach---which we refer to as `OLR-BT'---combines balanced truncation with the optimal low-rank approach for posterior approximation from \cite{Spantini2015Optimal}. PD-BT and OLR-BT share the same observability Gramian $\mbf{Q}_\infty$; compare \eqref{eq:Inf_Gramians} and \cite[equation (3.2)]{Qian2021Balancing}.

Key aspects in which PD-BT and OLR-BT differ are the role of the prior covariance $\PrCov$, the choice of input port, and the reachability Gramian.
PD-BT is motivated by balanced truncation of the prior-driven system \eqref{eq:PriorDrivenSys}, which matches the weighted output of the Bayesian smoothing problem \eqref{eq:BayesianSystem} \cite[Section 3.2.1]{2025PD-BT}. Note that in  \eqref{eq:PriorDrivenSys}, a square root $\mbf{L}_\pr$ of $\PrCov$ represents the input port, resulting in the reachability Gramian $\mbf{P}_\infty$ in \eqref{eq:Inf_Gramians}, which identifies easily reachable directions of the prior-driven system.
In contrast, OLR-BT assumes the existence of an input matrix $\mbf{B}$ such that $\mbf{A}\PrCov+\PrCov\mbf{A}^\top=-\mbf{BB}^\top$, so that the corresponding reachability Gramian is $\PrCov$; see \cite[Proposition 4.1]{Qian2021Balancing}.
Thus in OLR-BT, an easily reachable direction corresponds to a direction with high prior uncertainty, and these directions must be retained in the reduced model to be updated with measurements. The LTI system underlying OLR-BT is given by
\begin{align*}
    \dot{\mbf{x}}(t)  = \mbf{A}\mbf{x}(t)+\mbf{Bu}(t), \enspace \mbf{x}(0) = \mbf{0}; \quad
        \mbf{y} = \mbf{\Gamma}_\epsilon^{-1/2}\mbf{C}\mbf{x}(t),
\end{align*}
and has the impulse response
\begin{align}\label{eq:OLR impulse response}
    \mbf{h}_\mathrm{OLR}(t)= \mbf{\Gamma}_\epsilon^{-1/2}\mbf{C}e^{\mbf{A}t}\mbf{B}.   
\end{align}

Next, we show that we cannot use the strategy to prove error bounds for PD-BT to also prove error bounds for OLR-BT.
For PD-BT, an important step towards proving the error bounds in \Cref{thm:PD-BT bounds} was to exploit the fact that the rows of the difference between the exact and approximate root prior-preconditioned Hessians are exactly the differences of impulse responses corresponding to the full and the reduced systems; see \eqref{eq:PD-BT_reasoning}.
We now consider the analogous difference between the full and the reduced prior-preconditioned Hessian for OLR-BT. 
Let $\mbf{\widehat{G}}^\mathrm{OLR}$ be a reduced matrix of type \eqref{eq:PD-BT forward map}, but with system matrices $\mbf{C}^\mathrm{OLR}_r$ and $\mbf{A}^\mathrm{OLR}_r$ obtained from OLR-BT; see \cite[Section 3.2]{Qian2021Balancing}.
The analogue of \eqref{eq:PD-BT_reasoning} for OLR-BT then becomes
\begin{align}\label{eq:OLR-BT reasoning}
\begin{split}
        &\ObsCov^{-1/2}(\mbf{G-\widehat{G}^\mathrm{OLR}}) \mbf{L}_\pr\\
        &= \begin{bmatrix} \mbf{\Gamma}_\epsilon^{-1/2}\\ & \ddots \\ && \mbf{\Gamma}_\epsilon^{-1/2} \end{bmatrix}\cdot\Bigg(\begin{bmatrix}\mbf{C} e^{\mbf{A}t_1}\\ \vdots\\\mbf{C} e^{\mbf{A}t_n}\end{bmatrix}-\begin{bmatrix}\mbf{C}^\mathrm{OLR}_r e^{\mbf{A}^\mathrm{OLR}_rt_1}\\ \vdots\\\mbf{C}^\mathrm{OLR}_r e^{\mbf{A}^\mathrm{OLR}_rt_n}\end{bmatrix}(\mbf{V}^\mathrm{OLR}_r)^\top\Bigg) \cdot \begin{bmatrix}\mbf{L}_\pr \end{bmatrix}\\
        &= \begin{bmatrix} \mbf{\Gamma}_\epsilon^{-1/2}\mbf{C} e^{\mbf{A}t_1}\mbf{L}_\pr\\ \vdots\\  \mbf{\Gamma}_\epsilon^{-1/2}\mbf{C} e^{\mbf{A}t_n}\mbf{L}_\pr\end{bmatrix}-\begin{bmatrix} \mbf{\Gamma}_\epsilon^{-1/2}\mbf{C}_r e^{\mbf{A}_rt_1}\mbf{\Sigma}^{1/2}_r\mbf{Z}_r^\top\\ \vdots\\  \mbf{\Gamma}_\epsilon^{-1/2}\mbf{C}_r e^{\mbf{A}_rt_n}\mbf{\Sigma}^{1/2}_r\mbf{Z}_r^\top\end{bmatrix},    
\end{split}
\end{align}
where $\mbf{\Sigma}_r$ and $\mbf{Z}_r$ are obtained from a rank-$r$ truncated SVD of $\mbf{R}^\top \mbf{L}_\mathrm{pr}\approx\mbf{U}_r\mbf{\Sigma}_r\mbf{Z}_r^\top$, where $\mathbf{Q}_\infty =\mbf{R}\mbf{R}^\top$ and $\PrCov=\mbf{L}_\mathrm{pr}\mbf{L}_\mathrm{pr}^\top$.
For detailed definitions of the objects and explanations of the equations in \eqref{eq:OLR-BT reasoning}, see \cite[Section 2.3.2 and 3.2]{Qian2021Balancing}.
On the other hand, the input port $\mbf{B}$ and thus the impulse response $\mbf{h}_\mathrm{OLR}$ of OLR-BT given in \eqref{eq:OLR impulse response} do not appear in the prior-preconditioned Hessian, in contrast to the equation \eqref{eq:PD-BT_reasoning} for PD-BT.
As a result, for OLR-BT, one cannot combine a quantification of the MOR approximation quality in terms of the truncated HSVs with the general local Lipschitz bounds for linear Gaussian inverse problems from \Cref{theorem03}. Thus, one cannot bound the errors in the approximate posterior mean and covariance produced by OLR-BT in terms of the truncated HSVs. 

To summarize, despite their similarities as balanced truncation-based methods for solving the Bayesian smoothing problem, the OLR-BT and PD-BT methods are fundamentally different: OLR-BT aims to recover optimal low-rank posterior approximations, while PD-BT aims to recover system-theoretic invariants of the prior-driven system \eqref{eq:PriorDrivenSys}.

\paragraph{Interpolatory methods} The bounds in \Cref{thm:PD-BT bounds} exploit the strong connection \eqref{eq:PD-BT_reasoning} between the impulse response of the prior-driven system \eqref{eq:PriorDrivenSys} and the prior-preconditioned Hessian of the associated Bayesian smoothing problem \eqref{eq:BayesianSystem}, as well as the local Lipschitz bounds on the errors of the approximate posterior covariance and mean from \Cref{theorem03}.
It is of interest to determine whether similar connections hold for other system-theoretic MOR methods that can be applied to \eqref{eq:PriorDrivenSys}.

Interpolatory methods are a standard approach to match the impulse response of a full system with a reduced LTI system \cite{AntBeaGuc2020InterpolatoryMeths, Benner2017MORBook, GugercinIRKA2008}.
These methods interpolate the full-order transfer function, i.e. the Laplace transform of the impulse response, using a reduced-order transfer function constructed via projection. 
Rational interpolation methods typically operate in the frequency domain, where they interpolate the transfer function and possibly its derivative at a set of interpolants. 
However, we are not aware of any general $\mathcal{H}_2$-type bound---i.e., a bound on the $L_2$-norm of the impulse response error, similar to \eqref{eq:InhomImpulseErrorBound} or \eqref{eq:TLInhomImpulseErrorBound} for PD-(TL)BT---for these methods.
Only $\mathcal{H}_2$-optimality criteria are available, see e.g. \cite{BUNSEGERSTNER2010h2optimalMIMO,GugercinIRKA2008,VANDOOREN2008H2optimalMIMO}. 
Such a bound would be sufficient to obtain posterior error bounds using \Cref{theorem03}, as we did in \cref{sec:PD_BT_bound}. 

For the prior-driven system, the impulse input at $t = 0$ can be considered \cite{2025PD-BT}.
In the frequency domain, such an impulse input corresponds to a single interpolation point at infinity. 
For transfer function expansion at infinity, one refers to the moments as `Markov parameters', and to the resulting problem as `partial realization'; see \cite[Section 11.2]{Antoulas2005Book}. 
Whereas partial realization seems a promising MOR approach for the prior-driven system \eqref{eq:PriorDrivenSys}, we are not aware of a $\mathcal{H}_2$-type bound for this approach, and thus cannot analyze this approach using the proof strategy that we applied in \cref{sec:PD_BT_bound}. 
The development of interpolatory methods that provide good approximations of the posterior, and the analysis of the corresponding posterior approximation quality, constitute interesting directions for future research.

\section{Conclusion}\label{sec:Conclusion}
In \Cref{thm:PD-BT bounds} and \Cref{thm:PD-TLBT bounds}, we presented the first error bounds on the mean and covariance of a Gaussian posterior resulting from MOR techniques applied to Bayesian inference, namely for prior-driven balanced truncation and prior-driven, time-limited balanced truncation. 
These bounds leverage the relationship between the prior-preconditioned Hessian of the Bayesian inverse problem and the impulse response of the prior-driven system.
The bounds also relate the posterior approximation error to the truncated Hankel singular values, which are a classic system-theoretic measure.
This emphasizes the robust link between systems theory for the prior-driven system and Bayesian inference for general LTI forward models, including challenging cases such as unstable systems and singular prior covariances. 

The significance of our general local Lipschitz stability result in \Cref{theorem03} is not to provide tight control of the error in the posterior mean and covariance, but rather to emphasize that the error is of first order in the error of the square roots of the prior-preconditioned Hessian. This shows that the approximation error of the forward model is filtered through the square root of the prior covariance and the noise precision.
Our numerical results indicate that bounds accurately capture the trends of the error as a function of the rank of the reduced system.

Our work states the first local Lipschitz stability result for linear Gaussian inverse problems and perturbations of the forward model, and shows how this result can be applied to smoothing problems for LTI forward models \eqref{eq:BayesianSystem}. 
For future research, it would be of interest to extend this result to analyze other balanced-truncation-based MOR methods for the Bayesian smoothing problem that we did not consider in this work, e.g. \cite{Koe2022TLBTDA,koenig2022TLBT4DVAR,Qian2021Balancing}. Other interesting research directions include applying system-theoretic MOR methods, such as interpolatory methods \cite{AntBeaGuc2020InterpolatoryMeths, Benner2017MORBook, GugercinIRKA2008}, to the prior-driven system. It would also be of interest to apply the local Lipschitz stability results of \cite{CvetkovicLie2025,Sprungk2020LocalLipschitz} to MOR methods for possibly nonlinear or non-Gaussian inverse problems.

\section*{Funding}
This work was partially funded by the Deutsche Forschungsgemeinschaft (DFG) --- Project-ID 318763901 -- SFB1294.

\section*{Acknowledgments}
The authors thank Giuseppe Carere (Uni. Potsdam) for suggesting this collaboration, Elizabeth Qian (Georgia Tech.) for many fruitful discussions about her work and other MOR methods, and Melina Freitag (Uni. Potsdam) for her feedback on the manuscript. The authors thank Christopher Beattie (Virginia Tech.) for making us aware of an error in a previous version of this work.

\bibliographystyle{siamplain}

\appendix
\section{Auxiliary results}
\label{section_proofs_auxiliary_lemmas}

In this section, we collect auxiliary statements that we shall use in the proof of \Cref{theorem03}.
The first statement is an elementary decomposition.
 \begin{lemma} 
 \label{lemma02}
Let $m,n\in\N$. Let $\mbf A_i,\mbf{B}_i\in\R^{m\times n}$ for $i=1,2$, $\mbf M\in\R^{m\times n}$, and let $\mbf b\in\R^{n}$ and $\mbf v\in\R^{n}$.
 Then 
 \begin{align*}
&\mbf A_1\mbf B_1^\top \mbf M(\mbf b-\mbf B_1\mbf v)-\mbf A_2\mbf B_2^\top \mbf M(\mbf b-\mbf B_2\mbf v)
\\
=&(\mbf A_1-\mbf A_2) \mbf B_1^\top \mbf M(\mbf b-\mbf B_1\mbf v)+\mbf A_2 \mbf B_1^\top \mbf M (\mbf B_2-\mbf B_1)\mbf v
\\
&+\mbf A_2(\mbf B_1-\mbf B_2)^\top \mbf M(\mbf b-\mbf B_2\mbf v).
\end{align*}
 \end{lemma}
\begin{proof}[Proof of \Cref{lemma02}]
We have
\begin{align*}
&\mbf A_1\mbf B_1^\top \mbf M(\mbf b-\mbf B_1\mbf v)-\mbf A_2\mbf B_2^\top \mbf M(\mbf b-\mbf B_2\mbf v)
\\
=& \mbf A_1\mbf B_1^\top \mbf M(\mbf b-\mbf B_1\mbf v)-\mbf A_2 \mbf B_1^\top \mbf M(\mbf b-\mbf B_1\mbf v)+\mbf A_2\mbf B_1^\top \mbf M(\mbf b-\mbf B_1\mbf v)
\\
&-\mbf A_2\mbf B_2^\top \mbf M(\mbf b-\mbf B_2\mbf v)
\\
=&(\mbf A_1-\mbf A_2) \mbf B_1^\top \mbf M(\mbf b-\mbf B_1\mbf v)+\mbf A_2 (\mbf B_1^\top \mbf M (\mbf b-\mbf B_1\mbf v)-\mbf B_2^\top \mbf M(\mbf b-\mbf B_2\mbf v))
\\
=&(\mbf A_1-\mbf A_2) \mbf B_1^\top \mbf M(\mbf b-\mbf B_1\mbf v)+\mbf A_2 \bigr(\mbf B_1^\top \mbf M (\mbf b-\mbf B_1\mbf v)
\\
&-\mbf B_1^\top \mbf M(\mbf b-\mbf B_2\mbf v)+\mbf B_1^\top \mbf M(\mbf b-\mbf B_2\mbf v)-\mbf B_2^\top \mbf M(\mbf b-\mbf B_2\mbf v)\bigr)
\\
=&(\mbf A_1-\mbf A_2) \mbf B_1^\top \mbf M(\mbf b-\mbf B_1\mbf v)+\mbf A_2 \mbf B_1^\top \mbf M (\mbf B_2-\mbf B_1)\mbf v
\\
&+\mbf A_2(\mbf B_1-\mbf B_2)^\top \mbf M(\mbf b-\mbf B_2\mbf v).
\end{align*}
 \end{proof}
Recall that for a matrix $\mbf{A}$, $\norm{\mbf{A}}_p$ denotes the $p$-Schatten norm of $\mathbf{A}$; see \Cref{section_notation}.
\begin{lemma}
 \label{lemma05}
 Let $\mbf B_i\in\R^{m\times n}$, $i=1,2$. Then for $1\leq p\leq \infty$,
 \begin{align*}
  &\Norm{\mbf B_1^\top (\mbf I+ \mbf B_1\mbf B_1^\top)^{-1}\mbf B_1-\mbf B_2^\top (\mbf I+ \mbf B_2\mbf B_2^\top)^{-1}\mbf B_2}_p\leq D \norm{\mbf B_1-\mbf B_2}_p,
  \\
  D =& \Norm{ (\mbf I+ \mbf B_1\mbf B_1^\top)^{-1}\mbf B_1}_\infty+\Norm{\mbf B_2}_\infty\Norm{\mbf B_1}_\infty\left(\Norm{\mbf B_1}_\infty+\Norm{\mbf B_2}_\infty\right)
  \\
  &+\Norm{(\mbf I+ \mbf B_2\mbf B_2^\top)^{-1}\mbf B_2}_\infty.
 \end{align*}
\end{lemma}
For the proof of \Cref{lemma05}, we shall use the following theorem.
 \begin{theorem}[{\cite[Theorem 4.1]{Wedin1973}}]
  \label{lemma01}
  Let $m\in\N$ and let $\mbf{A}_1,\mbf{A}_2\in \R^{m\times m}$ be invertible.
  Then for any unitarily invariant norm $\norm{\cdot}_\bullet$,
  \begin{equation*}
   \norm{\mbf{A}_1^{-1}-\mbf{A}_2^{-1}}_\bullet \leq \norm{\mbf{A}_1^{-1}}_\infty \norm{\mbf{A}_2^{-1}}_\infty \norm{\mbf{A}_1-\mbf{A}_2}_\bullet
  \end{equation*}
 \end{theorem}
  Note that the spectral norm is denoted by $\norm{\cdot}_2$ in \cite{Wedin1973}, whereas we use $\norm{\cdot}_2$ to denote the 2-Schatten norm, i.e. the Frobenius or Hilbert--Schmidt norm.
 We shall also use the following properties of the $p$-Schatten norms for $1\leq p\leq\infty$: for every $\mbf A$ and $\mbf B$ such that $\mbf{AB}$ is defined,
\begin{equation}
 \label{eq:property_Frobenius_and_operator_norm}
 \norm{\mbf{AB}}_p\leq \norm{\mbf A}_\infty\norm{\mbf B}_p,\qquad \norm{\mbf A}_\infty \leq \norm{\mbf A}_p.
 \end{equation}
 Recall also that $\norm{\mbf{A}^\top}_p=\norm{\mbf{A}}_p$ for any admissible $\mbf{A}$ and $p$.
\begin{proof}[Proof of \Cref{lemma05}]
 By adding zero twice, we have
 \begin{align*}
  &\mbf B_1^\top (\mbf I+ \mbf B_1\mbf B_1^\top)^{-1}\mbf B_1-\mbf B_2^\top (\mbf I+ \mbf B_2\mbf B_2^\top)^{-1}\mbf B_2
  \\
  =& \mbf B_1^\top (\mbf I+ \mbf B_1\mbf B_1^\top)^{-1}\mbf B_1-\mbf B_2^\top (\mbf I+ \mbf B_1\mbf B_1^\top)^{-1}\mbf B_1+ \mbf B_2^\top (\mbf I+ \mbf B_1\mbf B_1^\top)^{-1}\mbf B_1
  \\
  &-\mbf B_2^\top (\mbf I+ \mbf B_2\mbf B_2^\top)^{-1}\mbf B_1+\mbf B_2^\top (\mbf I+ \mbf B_2\mbf B_2^\top)^{-1}\mbf B_1-\mbf B_2^\top (\mbf I+ \mbf B_2\mbf B_2^\top)^{-1}\mbf B_2.
  \end{align*}
  Thus, by applying the $\norm{\cdot}_p$ norm, the triangle inequality, and \eqref{eq:property_Frobenius_and_operator_norm}, we obtain
  \begin{align*}
   &\Norm{\mbf B_1^\top (\mbf I+ \mbf B_1\mbf B_1^\top)^{-1}\mbf B_1-\mbf B_2^\top (\mbf I+ \mbf B_2\mbf B_2^\top)^{-1}\mbf B_2}_p
   \\
   \leq & \Norm{(\mbf B_1-\mbf B_2)^\top (\mbf I+ \mbf B_1\mbf B_1^\top)^{-1}\mbf B_1}_p +\Norm{\mbf B_2^\top \left((\mbf I+ \mbf B_1\mbf B_1^\top)^{-1}-(\mbf I+ \mbf B_2\mbf B_2^\top)^{-1}\right)\mbf B_1}_p
   \\
   &+ \Norm{\mbf B_2^\top (\mbf I+ \mbf B_2\mbf B_2^\top)^{-1}\left(\mbf B_1-\mbf B_2\right)}_p
   \\
   \leq & \Norm{\mbf B_1-\mbf B_2}_p\Norm{ (\mbf I+ \mbf B_1\mbf B_1^\top)^{-1}\mbf B_1}_\infty
   \\
   &+\Norm{\mbf B_2}_\infty\Norm{(\mbf I+ \mbf B_1\mbf B_1^\top)^{-1}-(\mbf I+ \mbf B_2\mbf B_2^\top)^{-1}}_p\Norm{\mbf B_1}_\infty
   \\
   &+ \Norm{\mbf B_2^\top (\mbf I+ \mbf B_2\mbf B_2^\top)^{-1}}_\infty\Norm{\mbf B_1-\mbf B_2}_p.
  \end{align*}
Next, we bound the second term on the right-hand side of the last inequality.
Since the Schatten norms are unitarily invariant, we may apply \Cref{lemma01} with $\norm{\cdot}_\bullet\leftarrow \norm{\cdot}_p$ and $\mbf A_i\leftarrow \mbf I+ \mbf B_i\mbf B_i^\top$ for $i=1,2$:
\begin{align*}
& \Norm{(\mbf I+ \mbf B_1\mbf B_1^\top)^{-1}-(\mbf I+ \mbf B_2\mbf B_2^\top)^{-1}}_p
\\
\leq & \Norm{(\mbf I+ \mbf B_1\mbf B_1^\top)^{-1}}_\infty\Norm{(\mbf I+ \mbf B_2\mbf B_2^\top)^{-1}}_\infty\Norm{\mbf B_1\mbf B_1^\top-\mbf B_2\mbf B_2^\top}_p.
\end{align*}
By applying the triangle inequality and \eqref{eq:property_Frobenius_and_operator_norm},
\begin{align*}
 \Norm{\mbf B_1\mbf B_1^\top-\mbf B_2\mbf B_2^\top}_p=& \Norm{\mbf B_1\mbf B_1^\top-\mbf B_2\mbf B_1^\top+\mbf B_2\mbf B_1^\top-\mbf B_2\mbf B_2^\top}_p
 \\
 \leq & \Norm{\mbf B_1-\mbf B_2}_p\Norm{\mbf B_1}_\infty +\Norm{\mbf B_2}_\infty\Norm{\mbf B_1-\mbf B_2}_p
 \\
 =&\Norm{\mbf B_1-\mbf B_2}_p\left(\Norm{\mbf B_1}_\infty+\Norm{\mbf B_2}_\infty\right).
\end{align*}
Combining the preceding bounds completes the proof of \Cref{lemma05}.
\end{proof}

\section{Proof of posterior error bounds}
\label{section_proofs_local_lipschitz_bounds}

\errorBoundppHessianPerturbation*

\begin{proof}[Proof of \Cref{theorem03}]
 We first prove the bound on the difference between the posterior covariances.
By \Cref{assumption_main}, $\ObsCov^{-1/2}$ exists. Since $ \mbf{L}_\pr \mbf{L}_\pr^\top=\PrCov$,
 \begin{align}
  &\PrCov\fwdmap^\top(\ObsCov+\fwdmap \PrCov\fwdmap^\top)^{-1}\fwdmap\PrCov
  \nonumber\\
  =&\mbf{L}_\pr\mbf{L}_\pr^\top\fwdmap^\top\left(\ObsCov^{1/2}(\mbf I+\ObsCov^{-1/2}\fwdmap \mbf{L}_\pr\mbf{L}_\pr^\top \fwdmap^\top \ObsCov^{-1/2})\ObsCov^{1/2}\right)^{-1}\fwdmap\mbf{L}_\pr\mbf{L}_\pr^\top
  \nonumber\\
  =&\mbf{L}_\pr\mbf{L}_\pr^\top\fwdmap^\top\ObsCov^{-1/2} \left( \mbf I+\ObsCov^{-1/2}\fwdmap \mbf{L}_\pr\mbf{L}_\pr^\top \fwdmap^\top \ObsCov^{-1/2}\right)^{-1}\ObsCov^{-1/2}\fwdmap\mbf{L}_\pr\mbf{L}_\pr^\top
  \label{eq:11a_chol}
 \end{align}
 Replacing $\fwdmap$ with $\hatFwdmap$ in \eqref{eq:11a_chol} yields
\begin{align}
 &\PrCov\hatFwdmap^\top (\ObsCov+\hatFwdmap \PrCov\hatFwdmap^\top)^{-1} \hatFwdmap\PrCov
 \nonumber\\
 =&\mbf{L}_\pr\mbf{L}_\pr^\top\hatFwdmap^\top\ObsCov^{-1/2} \left( \mbf I+\ObsCov^{-1/2}\hatFwdmap \mbf{L}_\pr\mbf{L}_\pr^\top \hatFwdmap^\top \ObsCov^{-1/2}\right)^{-1}\ObsCov^{-1/2}\hatFwdmap\mbf{L}_\pr\mbf{L}_\pr^\top.
 \label{eq:11b_chol}
\end{align}
Thus
 \begin{align*}
&\Norm{\PosCov-\hatPosCov}_p
\\
=&\Norm{\PrCov \fwdmap^\top(\ObsCov+\fwdmap \PrCov \fwdmap^\top)^{-1}\fwdmap\PrCov -\PrCov \hatFwdmap^\top(\ObsCov+\hatFwdmap \PrCov \hatFwdmap^\top)^{-1}\hatFwdmap \PrCov}_p
\\
\leq &\Norm{\mbf{L}_\pr}_\infty^2 \left\Vert \mbf{L}_\pr^\top\fwdmap^\top\ObsCov^{-1/2} \left( \mbf I+\ObsCov^{-1/2}\fwdmap \mbf{L}_\pr\mbf{L}_\pr^\top \fwdmap^\top \ObsCov^{-1/2}\right)^{-1}\ObsCov^{-1/2}\fwdmap\mbf{L}_\pr\right. 
\\
&\left. \hskip10ex-\mbf{L}_\pr^\top\hatFwdmap^\top\ObsCov^{-1/2} \left( \mbf I+\ObsCov^{-1/2}\hatFwdmap \mbf{L}_\pr\mbf{L}_\pr^\top \hatFwdmap^\top \ObsCov^{-1/2}\right)^{-1}\ObsCov^{-1/2}\hatFwdmap\mbf{L}_\pr\right\Vert_p.
\end{align*}
 The equation follows by \eqref{eq:pos_covariance} and \eqref{eq:pos_covariance_approx}.
 The inequality follows by applying \eqref{eq:11a_chol} and \eqref{eq:11b_chol}, applying \eqref{eq:property_Frobenius_and_operator_norm} twice to extract $\norm{\mbf{L}_\pr}_\infty$ and $\norm{\mbf{L}_\pr^\top}_\infty$, and using that $\norm{\mbf{L}_\pr^\top}_p=\norm{\mbf{L}_\pr}_p$ for every $1\leq p\leq\infty$.
 By applying \Cref{lemma05} with the substitutions $\mbf{B}_1\leftarrow \ObsCov^{-1/2}\fwdmap \mbf{L}_\pr$ and $\mbf{B}_2\leftarrow \ObsCov^{-1/2} \hatFwdmap \mbf{L}_\pr$, we can bound the second term on the right-hand side of the inequality to obtain 
 \begin{equation*}
  \norm{\PosCov-\hatPosCov}_p\leq C \norm{\ObsCov^{-1/2} \fwdmap \mbf{L}_\pr-\ObsCov^{-1/2} \hatFwdmap \mbf{L}_\pr}_p,
 \end{equation*}
 where $C$ is obtained by making the above substitutions in the definition of the scalar $D$ from \Cref{lemma05}:
 \begin{align}
     C=& \Norm{\mbf{L}_{\pr}}_\infty^2 \biggr(\Norm{ (\mbf I+ \ObsCov^{-1/2}\fwdmap \PrCov\fwdmap^\top \ObsCov^{-1/2})^{-1}\ObsCov^{-1/2}\fwdmap \mbf{L}_\pr}_\infty
 \label{eq:error_covariance_for_perturbed_forward_model_by_ppH_constant} \\
  &+\Norm{\ObsCov^{-1/2}\hatFwdmap \mbf{L}_\pr}_\infty\Norm{\ObsCov^{-1/2}\fwdmap \mbf{L}_\pr}_\infty\left(\Norm{\ObsCov^{-1/2}\fwdmap \mbf{L}_\pr}_\infty+\Norm{\ObsCov^{-1/2}\hatFwdmap \mbf{L}_\pr}_\infty\right)
  \nonumber\\
  &+\Norm{(\mbf I+ \ObsCov^{-1/2}\hatFwdmap \PrCov \hatFwdmap^\top\ObsCov^{-1/2})^{-1}\ObsCov^{-1/2}\hatFwdmap \mbf{L}_\pr}_\infty\biggr).
\nonumber
 \end{align}
 This proves the bound \eqref{eq:error_covariance_for_perturbed_forward_model_by_ppH} on the $p$-Schatten error of $\hatPosCov$.
 
 To prove the bound on the $\ell_2$ error of $\hatPosMn$, we note that
 \begin{align}
 \PosMn(\data)-\hatPosMn(\data) =&  \PosCov \fwdmap^\top  \ObsCov ^{-1}(\data-\fwdmap \PrMn )- \hatPosCov\hatFwdmap^\top  \ObsCov ^{-1}(\data-\hatFwdmap\PrMn ) 
 \nonumber\\
 =&( \PosCov -\hatPosCov)\fwdmap^\top \ObsCov ^{-1}(\data-\fwdmap \PrMn )+\hatPosCov\fwdmap^\top \ObsCov ^{-1}(\hatFwdmap-\fwdmap)\PrMn  
 \label{eq:08_chol}\\
 &+\hatPosCov(\fwdmap-\hatFwdmap)^\top \ObsCov ^{-1}(\data-\hatFwdmap\PrMn ) 
\nonumber
\end{align}
where the first equation follows by \eqref{eq:pos_mean} and \eqref{eq:pos_mean_approx}, and the second equation follows by the application of \Cref{lemma02} with the substitutions $\mbf A_1\leftarrow  \PosCov $, $\mbf A_2\leftarrow \hatPosCov$, $\mbf B_1\leftarrow \fwdmap$, $\mbf B_2\leftarrow\hatFwdmap$, $\mbf M\leftarrow  \ObsCov ^{-1}$, $\mbf b\leftarrow \data$, and $\mbf v\leftarrow \PrMn $.

We now bound the second term on the right-hand side of \eqref{eq:08_chol}.
By the hypothesis that $\PrMn\in\range \mbf{L}_\pr$, there exists some $\mbf z$ such that $\PrMn=\mbf{L}_\pr\mbf z$.
Combining this fact with the generalized inverse property of the Moore--Penrose inverse $\mbf{L}_\pr^\dagger$, i.e. the property that $\mbf{L}_\pr\mbf{L}_\pr^\dagger \mbf{L}_\pr=\mbf{L}_\pr$,
\begin{align*}
 \hatPosCov\fwdmap^\top \ObsCov ^{-1}(\hatFwdmap-\fwdmap)\PrMn  =&\hatPosCov\fwdmap^\top \ObsCov ^{-1}(\hatFwdmap-\fwdmap)\mbf{L}_\pr\mbf{L}_\pr^\dagger \mbf{L}_\pr\mbf z
 \\
 =&\hatPosCov\fwdmap^\top \ObsCov ^{-1/2}\ObsCov ^{-1/2}(\hatFwdmap-\fwdmap)\mbf{L}_\pr \mbf{L}_\pr^\dagger \PrMn.
\end{align*}

Now consider the third term on the right-hand side of \eqref{eq:08_chol}.
By \eqref{eq:pos_covariance_approx}, the hypothesis $\mbf{L}_\pr\mbf{L}_\pr^\top=\PrCov$, and the generalized inverse property $\mbf{L}_\pr\mbf{L}_\pr^\dagger \mbf{L}_\pr=\mbf{L}_\pr$, 
\begin{align*}
 \hatPosCov= &\PrCov-\PrCov \hatFwdmap^\top(\ObsCov+\hatFwdmap \PrCov \hatFwdmap^\top)^{-1} \hatFwdmap\PrCov
 \\
 =&(\mbf I-\PrCov \hatFwdmap^\top(\ObsCov+\hatFwdmap \PrCov \hatFwdmap^\top)^{-1} \hatFwdmap)\mbf{L}_\pr\mbf{L}_\pr^\top
 \\
 =&(\mbf I-\PrCov \hatFwdmap^\top(\ObsCov+\hatFwdmap \PrCov \hatFwdmap^\top)^{-1} \hatFwdmap)\mbf{L}_\pr\mbf{L}_\pr^\top(\mbf{L}_\pr^\dagger)^\top \mbf{L}_\pr^\top =\hatPosCov (\mbf{L}_\pr^\dagger)^\top\mbf{L}_\pr^\top.
\end{align*}
Thus, the third term on the right-hand side of \eqref{eq:08_chol} satisfies
\begin{equation*}
 \hatPosCov(\fwdmap-\hatFwdmap)^\top \ObsCov ^{-1}(\data-\hatFwdmap\PrMn ) =\hatPosCov (\mbf{L}_\pr^\dagger)^\top\mbf{L}_\pr^\top(\fwdmap-\hatFwdmap)^\top \ObsCov ^{-1/2}\ObsCov^{-1/2}(\data-\hatFwdmap\PrMn ) .
\end{equation*}
By combining the preceding equations with \eqref{eq:08_chol}, taking the $\norm{\cdot}_2$ norm, applying the triangle inequality, and recalling that $\norm{\mbf{A}\mbf{x}}_2\leq\norm{\mbf{A}}_\infty\norm{\mbf{x}}_2$ for admissible $\mbf{A}$ and $\mbf{x}$,
\begin{align*}
 &\Norm{\PosMn(\data)-\hatPosMn(\data)}_2
 \\
 \leq & \Norm{( \PosCov -\hatPosCov)\fwdmap^\top  \ObsCov ^{-1}(\data-\fwdmap \PrMn )}_2
 +\Norm{\hatPosCov\fwdmap^\top \ObsCov ^{-1/2}\ObsCov ^{-1/2}(\hatFwdmap-\fwdmap)\mbf{L}_\pr \mbf{L}_\pr^\dagger \PrMn}_2
 \\
 &+\Norm{\hatPosCov (\mbf{L}_\pr^\dagger)^\top\mbf{L}_\pr^\top(\fwdmap-\hatFwdmap)^\top \ObsCov ^{-1/2}\ObsCov^{-1/2}(\data-\hatFwdmap\PrMn ) }_2
 \\
\leq &\Norm{ \PosCov -\hatPosCov}_\infty \Norm{\fwdmap^\top  \ObsCov ^{-1}(\data-\fwdmap \PrMn )}_2
\\
&+\Norm{\hatPosCov\fwdmap^\top \ObsCov ^{-1/2}}_\infty\Norm{\ObsCov^{-1/2}(\hatFwdmap-\fwdmap)\mbf{L}_\pr}_\infty\Norm{\mbf{L}_\pr^\dagger \PrMn}_2
 \\
 &+\Norm{\hatPosCov (\mbf{L}_\pr^\dagger)^\top }_\infty\Norm{\mbf{L}^\top_\pr(\fwdmap-\hatFwdmap)^\top \ObsCov ^{-1/2}}_\infty\Norm{\ObsCov^{-1/2}(\data-\hatFwdmap\PrMn )}_2.
 \end{align*}
We bound the first term and the sum of the second and third terms on the right-hand side of the last inequality separately.
For the first term, we have $\norm{ \PosCov -\hatPosCov}_\infty\leq C\Norm{\mbf{L}_{\pr}}_\infty^2\norm{\ObsCov ^{-1/2}(\fwdmap-\hatFwdmap) \mbf{L}_\pr}_\infty$ by  \eqref{eq:error_covariance_for_perturbed_forward_model_by_ppH}, with $C$ as in \eqref{eq:error_covariance_for_perturbed_forward_model_by_ppH_constant}.
Since $\norm{\mbf{A}}_\infty=\norm{\mbf{A}^\top}_\infty$ for any $\mbf{A}$, the sum of the second and third terms is equal to $D'\norm{\ObsCov ^{-1/2}(\fwdmap-\hatFwdmap) \mbf{L}_\pr}_\infty$, where
\begin{equation*}
 D'=\norm{\hatPosCov\fwdmap^\top \ObsCov ^{-1/2}}_\infty\norm{\mbf{L}_\pr^\dagger \PrMn}_2+\norm{\hatPosCov (\mbf{L}_\pr^\dagger)^\top}_\infty\norm{\ObsCov^{-1/2}(\data-\hatFwdmap\PrMn )}_2.
\end{equation*}
Combining these bounds yields \eqref{eq:error_mean_for_perturbed_forward_model_by_ppH}, with
\begin{align}
 &C'= C\Norm{\fwdmap^\top  \ObsCov ^{-1}(\data-\fwdmap \PrMn )}_2 +\Norm{\hatPosCov\fwdmap^\top \ObsCov ^{-1/2}}_\infty\Norm{\mbf{L}_\pr^\dagger \PrMn}_2
 \label{eq:error_mean_for_perturbed_forward_model_by_ppH_constant}
 \\
 &\hskip8ex +\Norm{\hatPosCov (\mbf{L}_\pr^\dagger)^\top }_\infty\Norm{\ObsCov^{-1/2}(\data-\hatFwdmap\PrMn )}_2,
 \nonumber
\end{align}
for $C$ as in \eqref{eq:error_covariance_for_perturbed_forward_model_by_ppH_constant}. 
This completes the proof of \Cref{theorem03}.
\end{proof}
\end{document}